\documentclass[12pt]{article}
\usepackage{amsfonts,amsmath,amsxtra}
\usepackage{amssymb}
\def\hybrid{\topmargin 0pt      \oddsidemargin 0pt
        \headheight 0pt \headsep 0pt
        \textwidth 16.5cm
        \textheight 23cm
        \hoffset=0.4cm
        \marginparwidth 0.0in
        \parskip 5pt plus 1pt   \jot = 1.5ex}
\catcode`\@=11
\def\marginnote#1{}

\newcount\hour
\newcount\minute
\newtoks\amorpm
\hour=\time\divide\hour by60
\minute=\time{\multiply\hour by60 \global\advance\minute by-\hour}
\edef\standardtime{{\ifnum\hour<12 \global\amorpm={am}%
        \else\global\amorpm={pm}\advance\hour by-12 \fi
        \ifnum\hour=0 \hour=12 \fi
      \number\hour:\ifnum\minute<10 0\fi\number\minute\the\amorpm}}
\edef\militarytime{\number\hour:\ifnum\minute<10 0\fi\number\minute}

\def\draftlabel#1{{\@bsphack\if@filesw {\let\thepage\relax
   \xdef\@gtempa{\write\@auxout{\string
      \newlabel{#1}{{\@currentlabel}{\thepage}}}}}\@gtempa
   \if@nobreak \ifvmode\nobreak\fi\fi\fi\@esphack}
        \gdef\@eqnlabel{#1}}
\def\@eqnlabel{}
\def\@vacuum{}
\def\draftmarginnote#1{\marginpar{\raggedright\scriptsize\tt#1}}

\def\draft{\oddsidemargin -0.1truein
        \def\@oddfoot{\sl preliminary draft \hfil
        \rm\thepage\hfil\sl\today\quad\militarytime}
        \let\@evenfoot\@oddfoot \overfullrule 3pt
        \let\label=\draftlabel
        \let\marginnote=\draftmarginnote
\def\@eqnnum{{\rm (\theequation)}
\rlap{\kern\marginparsep\tt\@eqnlabel}%
\global\let\@eqnlabel\@vacuum}  }

\newcommand{\RR}{{\mathbb{R}}}
\newcommand{\CC}{{\mathbb{C}}}

\newfont{\Bbbb}{msbm7 scaled 1\@ptsize00}
\newcommand{\zs}{\raise-1pt\hbox{$\mbox{\Bbbb Z}$}}

\font\teneufm=cmmib10 scaled 1\@ptsize00
\font\seveneufm=cmmib7 scaled 1\@ptsize00
\font\fiveeufm=cmmib5  
\def\bfit#1{{\textfont1=\teneufm\scriptfont1=\seveneufm
\scriptscriptfont1=\fiveeufm
\mathchoice{
\hbox{$\mathsurround=0pt\displaystyle#1$}}
{\mathsurround=0pt\hbox{$\textstyle#1$}}
{\hbox{$\mathsurround=0pt\scriptstyle#1$}}
{\hbox{$\mathsurround=0pt\scriptscriptstyle#1$}}}}
\font\sevenmsa=msam6 
\newfam\msafam
\textfont\msafam=\sevenmsa
\def\hexnumber@#1{\ifnum#1<10 \number#1\else
\ifnum#1=10 A\else\ifnum#1=11 B\else\ifnum#1=12 C\else
\ifnum#1=13 D\else\ifnum#1=14 E\else\ifnum#1=15 F\fi\fi\fi\fi\fi\fi\fi}
\def\msa@{\hexnumber@\msafam}
\def\llcorner{\delimiter"4\msa@78\msa@78 }
\def\lrcorner{\delimiter"5\msa@79\msa@79 }
\mathchardef\blacktriangleright="3\msa@49
\mathchardef\blacktriangleleft="3\msa@4A
\mathchardef\trianglerighteq="3\msa@44
\mathchardef\trianglelefteq="3\msa@45
\font\tenmsb=msbm10 scaled 1\@ptsize00
\newfam\msbfam
\textfont\msbfam=\tenmsb
\scriptfont\msbfam=\tenmsb
\def\msb@{\hexnumber@\msbfam}
\mathchardef\varkappa="0\msb@7B

\newdimen\linethick  \linethick=0.4pt
\newdimen\hboxitspace    \hboxitspace=5pt
\newdimen\vboxitspace    \vboxitspace=5pt

\def\fr#1{%
\beq\new
\vcenter{
\hrule height\linethick
           \hbox{\vrule width\linethick
                 \kern\hboxitspace
                 \vbox{\kern\vboxitspace
                       \hbox{$\begin{array}{c}\displaystyle#1
          \end{array}$}%
                       \kern\vboxitspace}%
                 \kern\hboxitspace
                 \vrule width\linethick}%
           \hrule height\linethick}%
\eeq}

\newdimen\Squaresize \Squaresize=14pt
\newdimen\Thickness \Thickness=0.5pt

\def\Square#1{\hbox{\vrule width \Thickness
   \vbox to \Squaresize{\hrule height \Thickness\vss
      \hbox to \Squaresize{\hss#1\hss}
   \vss\hrule height\Thickness}
\unskip\vrule width \Thickness}
\kern-\Thickness}

\def\Vsquare#1{\vbox{\Square{$#1$}}\kern-\Thickness}

\def\numberbysection{\@addtoreset{equation}{section}
        \def\theequation{\thesection.\arabic{equation}}}
\numberbysection

\renewcommand{\theequation}{\thesection.\arabic{equation}}
\newcommand{\l@qq}[2]{\addvspace{2em}
 \hbox to\textwidth{\hspace{1em}\bf #1 \dotfill #2}}

\newcounter{app}

\def\app{\setcounter{equation}{0}
\def\theequation{\Alph{app}.\arabic{equation}}\par
   \addvspace{10ex}
   \@afterindentfalse
  \secdef\@app\@dapp}
\newcommand\@app{\@startsection {app}{1}{-0.3ex}%
                             {-3.5ex \@plus -1ex \@minus -.2ex}%
                                   {2.3ex \@plus.2ex}%
                                   {\normalfont\Large\bf}}

\def\@dapp#1{%
{\parindent \z@ \raggedright \bf #1}\par\nobreak}
\def\l@app#1#2{\ifnum \c@tocdepth >\z@
    \addpenalty\@secpenalty
    \addvspace{1.0em \@plus\p@}%
    \setlength\@tempdima{1.5em}%
    \begingroup
      \parindent \z@ \rightskip \@pnumwidth
      \parfillskip -\@pnumwidth
      \leavevmode \bfseries
      \advance\leftskip\@tempdima
      \hskip -\leftskip
      #1\nobreak\hfil \nobreak\hb@xt@\@pnumwidth{\hss #2}\par
    \endgroup\fi}
\newcounter{sapp}[app]

\def\sapp{\def\theequation{\Alph{app}.\arabic{equation}}\par
   \@afterindentfalse
  \secdef\@sapp\@dsapp}
\newcommand\@sapp{\@startsection{sapp}{2}{\z@}%
                           {-3.25ex\@plus -1ex \@minus -.2ex}%
                           {1.5ex \@plus .2ex}%
                              {\normalfont\large\bfseries}}

\def\@dsapp#1{%
{\parindent \z@ \raggedright  \bf #1}\par\nobreak}

\newcommand{\l@sapp}{\@dottedtocline{2}{1.4em}{2.5em}}

\def\titlepage{\@restonecolfalse\if@twocolumn\@restonecoltrue\onecolumn
     \else \newpage \fi \thispagestyle{empty}\c@page\z@
        \def\thefootnote{\fnsymbol{footnote}} }

\def\endtitlepage{\if@restonecol\twocolumn \else  \fi
        \def\thefootnote{\arabic{footnote}}
        \setcounter{footnote}{0}}  
\relax

\hybrid

\newtoks\@stequation

\def\subequations{\refstepcounter{equation}%
  \edef\@savedequation{\the\c@equation}%
  \@stequation=\expandafter{\theequation}
  \edef\@savedtheequation{\the\@stequation}
  \edef\oldtheequation{\theequation}%
  \setcounter{equation}{0}%
  \def\theequation{\oldtheequation\alph{equation}}}

\def\endsubequations{%
  \setcounter{equation}{\@savedequation}%
  \@stequation=\expandafter{\@savedtheequation}%
  \edef\theequation{\the\@stequation}%
  \global\@ignoretrue}

\makeatletter
\newdimen\normalarrayskip            
\newdimen\minarrayskip               
\normalarrayskip\baselineskip
\minarrayskip\jot
\newif\ifold             \oldtrue            \def\new{\oldfalse}
\def\arraymode{\ifold\relax\else\displaystyle\fi}
\def\eqnumphantom{\phantom{(\theequation)}} 
\def\@arrayskip{\ifold\baselineskip\z@\lineskip\z@
     \else
     \baselineskip\minarrayskip\lineskip1\baselineskip\fi}

\def\@arrayclassz{\ifcase \@lastchclass \@acolampacol \or
\@ampacol \or \or \or \@addamp \or
   \@acolampacol \or \@firstampfalse \@acol \fi
\edef\@preamble{\@preamble
  \ifcase \@chnum
     \hfil$\relax\arraymode\@sharp$\hfil
     \or $\relax\arraymode\@sharp$\hfil
     \or \hfil$\relax\arraymode\@sharp$\fi}}

\def\@array[#1]#2{\setbox\@arstrutbox=\hbox{\vrule
     height\arraystretch \ht\strutbox
     depth\arraystretch \dp\strutbox
width\z@}\@mkpream{#2}\edef\@preamble{\halign \noexpand\@halignto
\bgroup \tabskip\z@ \@arstrut \@preamble \tabskip\z@ \cr}%
\let\@startpbox\@@startpbox \let\@endpbox\@@endpbox
  \if #1t\vtop \else \if#1b\vbox \else \vcenter \fi\fi
  \bgroup \let\par\relax
  \let\@sharp##\let\protect\relax
  \@arrayskip\@preamble}
\def\eqnarray{\stepcounter{equation}%
              \let\@currentlabel=\theequation
              \global\@eqnswtrue
              \global\@eqcnt\z@
              \tabskip\@centering              
              \let\\=\@eqncr
              $$%
            \halign to \displaywidth  \bgroup
             \eqnumphantom \@eqnsel
      \hskip\@centering                               
    $\displaystyle  \tabskip\z@ {##}$%
    &\global\@eqcnt\@ne \hskip 2\arraycolsep
         $ \displaystyle  \arraymode{##}$\hfil
    &\global\@eqcnt\tw@ \hskip 2\arraycolsep
         $\displaystyle\tabskip\z@{##}$\hfil
         \tabskip\@centering
    &{##}\tabskip\z@\cr}
\makeatother
\newtheorem{te}{Theorem}[section]
\newtheorem{de}{Definition}[section]
\newtheorem{prop}{Proposition}[section]           
\newtheorem{cor}{Corollary}[section]
\newtheorem{lem}{Lemma}[section]
\newtheorem{ex}{Example}[section]
\newtheorem{rem}{Remark}[section]

\newcommand{\beq}[1]{\begin{equation}\label{#1}}
\newcommand\eeq {\end{equation}}
\newcommand\bqa {\begin{eqnarray}}
\newcommand\eqa {\end{eqnarray}}
\def\be{\begin{eqnarray}\new\begin{array}{cc}}
\def\ee{\end{array}\end{eqnarray}}

\def\beq{\begin{equation}}
\def\eeq{\end{equation}}
\def\bse{\begin{subequations}}                
\def\ese{\end{subequations}}
\def\square{\hfill{\vrule height6pt width6pt            
depth1pt} \break \vspace{.01cm}}                        
\def\d{\partial}
\def\stack#1#2{\raise0.7pt\hbox{$\mathrel{\mathop{#2}\limits^{#1}}$}}

\def\bgamma{{\bfit\gamma}}

\def\brho{{\bfit\rho}}

\def\la{\lambda}

\def\<{\langle}
\def\>{\rangle}
\def\ov{\overline}

\def\frak{\mathfrak}

\def\ca{{\cal A}}
\def\cb{{\cal B}}
\def\cc{{\cal C}}

\def\N{{\scriptscriptstyle N}}
\def\ts#1#2{{\textstyle\frac{#1}{#2}}}

\def\ov{\overline}

\def\g{\gamma}

\begin{document}
\thispagestyle{empty}
\hfill {\normalsize ITEP-TH-20/02}

\begin{center}

\phantom.
\bigskip\bigskip
{\Large\bf  Representation Theory and the Quantum Inverse

\vspace{0.3cm}
Scattering Method:

\vspace{0.3cm}
The Open Toda Chain and the Hyperbolic Sutherland Model
}

\vspace{1cm}

\bigskip \bigskip
{\large A. Gerasimov\footnote{E-mail: gerasimov@vitep5.itep.ru},
 S. Kharchev\footnote{E-mail: kharchev@gate.itep.ru}, D.
Lebedev\footnote{E-mail: lebedev@gate.itep.ru}}\\ \medskip
{\it Max-Planck-Institut f\"ur Mathematik, Vivatsgasse 7,
D-53111 Bonn, Germany\\
and\\
Institute for Theoretical \& Experimental Physics, 117259, Moscow,
Russia}\\
\end{center}

\vspace{1.5cm}

\begin{abstract}
\noindent
Using the representation theory of $\frak{gl}(N,\RR)$, we express the
wave function of the $GL(N,\RR)$ Toda chain, which two of us recently
obtained by the Quantum Inverse Scattering Method, in terms of multiple
integrals. The main tool  is our generalization of the Gelfand-Zetlin
method to the case of infinite-dimensional representations of
$\frak{gl}(N,\RR)$.
The interpretation of this generalized construction in terms of the
coadjoint orbits is given and the connection with the Yangian
$Y(\frak{gl}(N))$ is discussed. We also give the hyperbolic Sutherland model
eigenfunctions expressed in terms of integrals in the Gelfand-Zetlin
representation. Using the example of the open Toda chain, we discuss
the connection between the Quantum Inverse Scattering Method and
Representation Theory.
\end{abstract}

\clearpage \newpage
\setcounter{footnote}0


\tableofcontents
\normalsize
\section{Introduction}

In the late seventies two basic approaches to quantum integrable
systems were formulated. Olshanetsky and Perelomov discovered a
direct relation between the representation theory of non-compact
semi-simple Lie groups and Sutherland-Calogero-Moser (SCM)
systems \cite{OP1}-\cite{OP3}. The methods of Representation Theory
provide a complete solution to the quantum SCM system.
Soon after, Kostant and Kazhdan realized that more general
classes of integrable systems, including open Toda chains, can be
solved by the same method \cite{Ko1}, \cite{Ko2}.
The wave functions of the SCM systems and of Toda
chains appear as zonal spherical functions \cite{Gelf}, \cite{HC},
and Whittaker functions \cite{Jac}-\cite{GKM}, respectively. There is
a generalization of the above method to a still more general class of
integrable systems, formally  corresponding to affine Lie algebras and
affine quantum groups (see \cite{EFK}, \cite{E} and references therein).
However, up to now this approach has not produced explicit
representations of the wave functions in terms of integrals because
of analytic difficulties in the corresponding representation theory.

At the same time, another approach to quantum integrability
emerged. In the works of Faddeev, Sklyanin, Kulish and Takhtajan
\cite{FS1}-\cite{FS4}
\footnote{
See the summary of the
early period of the development in \cite{F1}.
},
the Quantum Inverse Scattering Method (QISM) was formulated and
successively applied to  numerous quantum integrable systems \cite{KS}.
A related recent development is connected with the Separation of
Variables (SoV) method \cite{Skl2}, \cite{Skl1}. Nevertheless,  the
interrelations of QISM with Representation Theory have  remained unclear.
The main objective of this paper is to outline an approach towards
establishing such a connection.

Our starting point was a new method of deriving the wave function
for the $GL(N,\RR)$ Toda chain given in \cite{KL2}. This had no obvious
relationship with the Representation Theory approach, but instead
QISM was used to construct the wave function.
The understanding of this construction in terms of  Representation Theory,
and the comparison of these two basic approaches to quantum
integrability for this simple integrable system, should provide additional
insight into the nature of QISM.
In this paper we give an interpretation of the integral expression
of the wave function of the open Toda chain, obtained in  \cite{KL2},
in terms of the Representation Theory approach to quantum
integrability. The basic tool is a generalization of the Gelfand-Zetlin
method \cite{GelZ}, \cite{GelG} to the case of infinite-dimensional
representations of $\frak{gl}(N,\RR)$.
This gives a construction of a representation of
${\cal U}(\frak{gl}(N,\RR))$ in terms of difference operators
\footnote{Incidentally,
a version of this construction for the affine algebra
$\widehat{\frak{sl}(2)}_k$ has been
discovered previously in \cite{GMM}. The generalization to the case
$\widehat{\frak{sl}(N)}_k$ will be given elsewhere.
},
which one can use for finding explicit expressions for Whittaker
vectors. Thus, we reproduce the  integral formula
for Whittaker functions obtained previously in \cite{KL2} by  QISM.

Comparing the methods of \cite{KL2} with those proposed in
this paper reveals a deep connection between the SoV
\cite{Skl2} and Gelfand-Zetlin type representations. In
\cite{FM}, \cite{Skl2} the separated variables for
the classical Toda chain were found. They could be successfully
quantized to give an explicit solution to the corresponding quantum
theory. In this paper we use the generalization of the Gelfand-Zetlin
representation for $\frak{gl}(N,\RR)$, which gives rise to a much larger
set of separated variables ($\frac{1}{2}N(N-1)$ instead of $(N-1)$).
The idea of introducing the enlarged set of separated variables can be
traced back to the work of Gelfand and Ki\-rillov on the  field
of fractions of the universal enveloping algebra \cite{GK1}.
On the other hand, the representation which we have constructed appears
to be closely connected with the representation theory of the Yangian
$Y(\frak{gl}(N))$. The enlargement of the set of variables
for solving the quantum problem is a crucial part of the Representation
Theory approach to quantum integrability. One may hope that understanding
the correct notion of separated variables in Representation Theory
will turn out to be fruitful for both the theory of Quantum Integrable
Systems and Representation Theory.

The plan of the paper is as follows. In section 2 we construct the main
tool we will be using throughout the paper:
we give a direct generalization of the Gelfand-Zetlin method
to the case of infinite-dimensional representations of
${\cal U}(\frak{gl}(N))$. The crucial point is to use continuous
unconstrained variables instead of integer-valued
Gelfand-Zetlin parameters. Another type of analytical continuation was
considered in \cite{GelG} (see also \cite{LP}).
The main results of this section are the analytic construction of a
representation of ${\cal U}(\frak{gl}(N))$ and explicit description of
mutually dual Whittaker modules.

In section 3 we give the interpretation of the representation described
in section 2 in terms of the Kirillov-Kostant  orbit method \cite{Kir}.
We give an explicit parameterization of the finite cover of the open part
of the coadjoint orbit of $GL(N,\RR)$ with stabilizer
$H_\N=GL(1,\RR)^{\otimes N}$. We show that this parameterization yields
Darboux coordinates with respect to the canonical Kirillov-Kostant
symplectic structure on the orbit. The above construction is a version of
a construction by Alekseev, Faddeev and Shatashvili \cite{AFS}, \cite{S}
adopted for non-compact groups. Our parameterization yields  the classical
counterparts for the generators of $\frak{gl}(N)$ given in section 2.
A subtle point of the orbit method interpretation is that we use only
the finite cover of the open part of the orbit. This leads to a
representation of $\frak{gl}(N,\RR)$ which cannot be integrated to
a representation of the  group. However, it is obvious that the action
of the Cartan subalgebra integrates to the action of a Cartan torus which
is sufficient for constructing solutions of the open Toda chain.
We end this section by  giving contour integral formulas for the
generators of $\frak{gl}(N,\RR)$.

In section 4 we relate the Gelfand-Zetlin parameterization to the
representation theory of the Yangian \cite{Dr2}. (Closely connected ideas
for the classical Gelfand-Zetlin construction were developed in
\cite{Cher1}-\cite{NT1}). In this section we derive explicit formulas
in the Gelfand-Zetlin parameterization for the images of the Drinfeld
generators $\ca_n(\la),\cb_n(\la),\cc_n(\la)$ under the natural
homomorphism: $Y(\frak{gl}(N))\rightarrow U(\frak{gl}(N))$.
This result will be used in section 6 to establish the connection
of the Gelfand-Zetlin parameterization with the QISM-based
approach to the integrability of the quantum open Toda chain.

In section 5 we prove two of the main results of this paper, Theorems 5.1
and 5.2, by using the methods developed in section 2.
We apply the representation constructed in section 2 to find the explicit
solutions of the quantum open Toda chain and the quantum hyperbolic
Sutherland model. For the Toda chain, we reproduce the results
of \cite{KL2}. This bridges the gap between the Representation Theory
and QISM approaches for these theories. The corresponding integral
representation for the wave functions of the quantum hyperbolic
Sutherland model appears to be new.

In the last section, we compare the QISM-based approach to the quantum
integrability  of the open Toda chain to the Representation Theory approach.

{\bf Acknowledgments:}
We would like to thank the Max Planck Institut f\"ur Mathematik
(Bonn, Germany), where the paper was finished, for financial support and
its stimulating atmosphere.

The authors are grateful to  V. Gorbunov, A. Levin, Y. Neretin,
M. Olshanetsky, and  M. Se\-me\-nov-Tian-Shansky for useful discussions.

The research was partly supported by grants RFBR 01-01-00548
(A. Gerasi\-mov); INTAS OPEN 00-00055; RFBR 00-02-16477, (S. Kharchev);
INTAS OPEN 00-00055; RFBR 01-01-00548 (D. Lebedev) and by
grant 00-15-96557 for Support of Scientific Schools.

\section{Whittaker modules in the Gelfand-Zetlin representation}
\subsection{ Generalization of the  Gelfand-Zetlin
construction to   infinite-\-di\-men\-sio\-nal
representations}

We start by recalling the original Gelfand-Zetlin construction
for the canonical basis in the space of an irreducible finite-dimensional
representation of $GL(N,\CC)$ \cite{GelZ}, \cite{GelG}.
It is  based  on the two  well-known  facts.

Any irreducible finite-dimensional representation is uniquely
determined by an $n$-dimen\-sio\-nal vector
$\bfit{m}_\N=(m_{\N 1},\ldots,m_{\N\N})$
with  integer entries such that
$m_{\N1}\geq\ldots\geq m_{\N\N}$. The vector $\bfit{m}_\N $ is
called the {\it label} or the highest weight of the representation.

Consider $GL(N-1,\CC)$ as the subgroup of $GL(N,\CC)$ embedded in
standard way. Then, each of the irreducible representations
of the subgroup  $GL(N-1,\CC)$ enters into the decomposition of a
finite-dimensional irreducible representation of  $GL(N,\CC)$
with multiplicity not greater than one. If we denote by
$\bfit{m}_{\N-1}$ the labels of the representations involved in the
decomposition, then their  entries  satisfy the  additional
conditions:
$m_{\N j}\geq m_{\N-1,j}\geq m_{\N,j+1}$, for all $j=1,\ldots,N-1$.

If we continue the decomposition step by step, we obtain a
decomposition of an irreducible finite-dimensional representation of
$GL(N,\CC)$ into a direct sum of one-dimensional subspaces. Fixing
an element in each one-dimensional summand  yields a Gelfand-Zetlin basis.
The vectors in this basis are labeled by triangular arrays, consisting
of all $\bfit{m}_n$, $n=1,\ldots,N$, with entries satisfying the above
mentioned constraints.

Now, we wish to generalize the Gelfand-Zetlin (GZ) construction to
infinite-dimen\-sio\-nal representations of the universal enveloping
algebra ${\cal U}(\frak{gl}(N))$.

We start with the following generalization of the Gelfand-Zetlin
representation of $\frak{gl}(N)$. Namely, let
$(\bgamma_1\,\ldots\,\bgamma_\N)$ be a triangular array consisting
of $\frac{1}{2}N(N+1)$ variables
$\bgamma_n=(\gamma_{n1},\ldots,\gamma_{nn})\in\CC^n;\,n=1,\ldots,N$.
\begin{prop}
Let $M$ be the space of meromorphic functions in $\frac{1}{2}N(N-1)$
variables $\bgamma_1,\ldots,\bgamma_{\N-1}$. Then the operators
\footnote{
Actually, we construct the family of representations depending on
the auxiliary parameter $\hbar$. One could rescale the variables
$\gamma_{nj}\to\hbar\gamma_{nj}$ to get rid of this parameter.
However, we will retain $\hbar$ in the formulas to make the connection
with standard notations in the theory of quantum integrable systems more
explicit.
}
\bse\label{m1}
\be\label{m1a}
\hspace{2cm}
E_{nn}=\frac{1}{i\hbar}\Big(\sum_{j=1}^n\gamma_{nj}-
\sum_{j=1}^{n-1}\gamma_{n-1,j}\Big);\hspace{2cm}(n=1,\ldots,N),
\ee
\be\label{m1b}
E_{n,n+1}=-\frac{1}{i\hbar}\sum_{j=1}^n\frac{\prod\limits_{r=1}^{n+1}
(\gamma_{nj}-\gamma_{n+1,r}-\frac{i\hbar}{2})}
{\prod\limits_{s\neq j}(\gamma_{nj}-\gamma_{ns})}\,
e^{-i\hbar\d_{\gamma_{nj}}};\hspace{0.5cm}(n=1,\ldots,N-1),
\hspace{-2cm}
\ee
\be\label{m1c}
E_{n+1,n}=\frac{1}{i\hbar}\sum_{j=1}^n\frac{\prod\limits_{r=1}^{n-1}
(\gamma_{nj}-\gamma_{n-1,r}+\frac{i\hbar}{2})}
{\prod\limits_{s\neq j}(\gamma_{nj}-\gamma_{ns})}\,
e^{i\hbar\d_{\gamma_{nj}}};\hspace{0.5cm}(n=1,\ldots,N-1)
\hspace{-2cm}
\ee
\ese
form  a  representation of $\frak{gl}(N)$ in $M$.
\end{prop}
{\bf Proof.} It is sufficient to show that the generators (\ref{m1})
satisfy the standard commutation relations
\be\label{s0}
[E_{nn},E_{m,m+1}]=(\delta_{nm}-\delta_{n,m+1})E_{m,m+1},\\

[E_{nn},E_{m+1,m}]=-(\delta_{nm}-\delta_{n,m+1})E_{m+1,m},\\

[E_{n,n+1},E_{m+1,m}]=(E_{nn}-E_{n+1,n+1})\delta_{nm}
\ee
and obey the Serre relations
\be\label{s1}
[E_{n,n+1},[E_{n,n+1},E_{n+1,n+2}]]\,=\,0,\\

[E_{n+1,n+2},[E_{n+1,n},E_{n+1,n+2}]]\,=\,0,\\

[E_{n+1,n},[E_{n+1,n},E_{n+2,n+1}]]\,=\,0,\\

[E_{n+2,n+1},[E_{n+1,n},E_{n+2,n+1}]]\,=\,0.
\ee
These follow by direct verification.
By defining the composite root generators recursively as
$E_{jk}=[E_{jm},E_{mk}]$ for $j< m<k$ and $j> m>k$, one obtains the
complete set of commutation relations
\be\label{alg}
[E_{jk},E_{lm}]=E_{jm}\delta_{lk}-E_{lk}\delta_{jm}.
\ee
\square

The representation described above is parameterized by the label
$\bgamma_\N=(\gamma_{\N1},\ldots,\gamma_{\N\N})$ and extends to
a representation of the universal enveloping algebra
${\cal U}(\frak{gl}(N))$. Obviously, there is a natural
representation of any subalgebra
${\cal U}(\frak{gl}(n))\subset{\cal U}(\frak{gl}(N)),\,
n=1,\ldots,N-1$ on the same space $M$ with the label
$\bgamma_n=(\gamma_{n1},\ldots,\gamma_{nn})$.
Let ${\cal Z}(\frak{gl}(n))$ be the centre of ${\cal U}(\frak{gl}(n))$.
\begin{de}
We say  that a  ${\cal U}(\frak{gl}(n))$-module $V$ admits an
infinitesimal character $\xi$ if there is a homomorphism
$\xi: {\cal Z}(\frak{gl}(n))\rightarrow \CC$ such that
$zv=\xi(z)v$ for all $z\!\in\!{\cal Z}(\frak{gl}(n)), v\in V.$
\end{de}
It is possible to show that the ${\cal U}(\frak{gl}(N))$-module $M$
defined above admits an infinitesimal character.
Actually, a more general statement holds:
\begin{prop}\label{lic}
Each  central element of ${\cal U}(\frak{gl}(n))$ acts on $M$  via
multiplication by a symmetric polynomial in the variables $\g_{nj}$.
\end{prop}
{\bf Proof.} Let us use the induction over $n$. For $n=1$ the
statement is obvious. Notice that the operators
$u\in{\cal U}(\frak{gl}(m))$ acting on $M$
are invariant under permutations in each subset of variables
$\!\bgamma_1,\ldots,\bgamma_m\!$.
Now suppose that the action of the centres
$\!{\cal Z}(\frak{gl}(m)),\, m=1,\ldots,n-1\!$ are given by the
functionally independent symmetric polynomials of
$\bgamma_m,\,m=1,\ldots,n-1$. Let us prove that the centres of
${\cal Z}(\frak{gl}(n))$ act on $M$ analogously.
The elements of the centre ${\cal Z}(\frak{gl}(n))$ commute with
${\cal Z}(\frak{gl}(m))$ and thus do not depend on the shifts over
$\bgamma_m$. On the other hand, the elements of ${\cal Z}(\frak{gl}(n))$
commute with all root generators $E_{ij},\,i,j=1,\ldots,n$ and, hence,
are $i\hbar$- periodic with respect to any variable
$\gamma_{mj},\,m=1,\ldots,n-1$. The polynomial generators of the
centre of ${\cal U}(\frak{gl}(n))$ may act only as rational
functions in this representation and thus we conclude that they do
not depend on the variables $\gamma_{mj},\,m=1,\ldots,n-1$. Taking
into account that the generators of ${\cal U}(\frak{gl}(n))$
depend polynomially on $\bgamma_n$ we conclude that the generators
of the centre ${\cal Z}(\frak{gl}(n))$ act on $M$ via  symmetric
polynomials in $\bgamma_n$. The functional independence may be directly
verified on a particular subset of $M$.\square

Hence, the elements $z\in{\cal Z}(\frak{gl}(n))$ act on the
${\cal U}(\frak{gl}(n))$-module $M$ as  operators of scalar multiplication,
and the ${\cal U}(\frak{gl}(n))$-module $M,$ $n=1,\ldots,N$, admits
an infinitesimal character.

In the next subsection we calculate the explicit action of the central
elements on $M$ using the notion of Whittaker vectors.

\subsection{Whittaker modules}

Let us now  give a construction for  Whittaker modules
using the representation of ${\cal U}(\frak{gl}(N))$ described above.
We first recall some facts from \cite{Ko2}.

Let $n_+$ and $n_-$ be the subalgebras of $\frak{gl}(N)$
generated,  respectively,  by positive and negative root generators.
The homomorphisms (characters) $\!\chi_+\!:n_+\rightarrow\CC\!$,
$\!\!\chi_-\!:n_-\!\!\rightarrow\CC$
are uniquely determined by their values on the simple root generators,
and are called {\it non-singular} if
the (complex) numbers $\chi_+(E_{n,n+1})$ and
$\chi_-(E_{n+1,n})$ are non-zero for all $n=1,\ldots,N-1$.

Let $V$ be any ${\cal U}={\cal U}(\frak{gl}(N))$-module. Denote  the
action of $u\in{\cal U}$ on $v\in V$ by $uv$.
A vector $w\in V$ is called a Whittaker vector with respect to the
character $\chi_+$ if
\be\label{ewt}
\hspace{2cm}
E_{n,n+1}w=\chi_+(E_{n,n+1})w\,,
\hspace{1cm}(n=1,\ldots,N-1),
\ee
and an element $w'\in V'$ is called a Whittaker vector
with respect to the character $\chi_-$  if
\be\label{fwt}
\hspace{2cm} E_{n+1,n}w'=\chi_-(E_{n+1,n})w'\,,
\hspace{1cm} (n=1,\ldots,N-1).
\ee

A Whittaker vector is {\it cyclic} for $V$ if ${\cal U}w=V$,
and a ${\cal U}$-module is a {\it Whittaker module}
if it contains a cyclic Whittaker vector. The
${\cal U}$-modules $V$ and $V'$ are called {\it dual}
if there exists a non-degenerate pairing $\<.\,,.\>: V'\times V\to\CC$
such that
$\langle X v',v\rangle = -\langle v',Xv\rangle$ for all $v\in V,$
$v'\in V'$ and  $X\in\frak{gl}(N)$.

We proceed with explicit formulas for Whittaker vectors
corresponding to the representation given by (\ref{m1}).
\begin{prop}
The equations
\be\label{fww'}
E_{n+1,n}w_\N' =-i\hbar^{-1}w_\N' ,
\ee
\be\label{fww''}
E_{n,n+1}w_\N=-i\hbar^{-1}w_\N
\ee
for all $n=1,\ldots,N-1$, admit the solutions
\be\label{wv}
w_\N'=1 ,\\
w_\N=e^{-\frac{\pi}{\hbar}\sum\limits_{n=1}^{N-1}(n-1)
\sum\limits_{j=1}^n
\gamma_{nj}}\prod_{n=1}^{N-1} s_n(\bgamma_n,\bgamma_{n+1}),
\ee
where
\be\label{wv'}
s_{n}(\bgamma_n,\bgamma_{n+1})=\prod_{k=1}^n\prod_{m=1}^{n+1}
\hbar^{\frac{\gamma_{nk}-\gamma_{n+1,m}}{i\hbar}+\frac{1}{2}}\;
\Gamma\Big(\frac{\gamma_{nk}-\gamma_{n+1,m}}{i\hbar}+\frac{1}{2}\Big).
\ee
\end{prop}
{\bf Proof.} The following equality holds:
\be\label{id1}
E_{n,n+1}w_\N=-i\hbar^{-1}w_\N\sum_{j=1}^n
\frac{\prod\limits_{r=1}^{n-1}
(\gamma_{nj}-\gamma_{n-1,r}-\frac{i\hbar}{2})}
{\prod\limits_{s\neq j}(\gamma_{nj}-\gamma_{ns})}\,.
\ee
Using the identity
\be\label{id2}
\sum_{j=1}^n\frac{x^m_j}{\prod\limits_{s\neq j}(x_j-x_s)}=
\sum_{\stackrel{\scriptstyle k_r\geq 0}{k_1+\ldots+k_n=m-n+1}}
x_1^{k_1}\ldots x_n^{k_n}
\ee
one arrives at (\ref{fww''}).  Equation (\ref{fww'}) is
proved similarly.
\square

Evidently, the solutions (\ref{wv}), (\ref{wv'}) are not unique.
The set of  Whittaker vectors is closed under multiplication
by an arbitrary $i\hbar$-periodic function in the  variables $\gamma_{nj}$.
Hence, there are infinitely many invariant subspaces in $M$ corresponding
to infinitely many Whittaker vectors.

To construct irreducible submodules, let us introduce the Whittaker
modules $W$ and $W',$ generated cyclically by the Whittaker vectors
$w_\N$ and $w_\N',$ respectively.
\begin{te}
Let $\bfit m_n=(m_{n1},\ldots,m_{nn})$ be the set of non-negative integers.
The Whittaker module $W={\cal U}w_\N$ is spanned by the elements
\be\label{bas1}
w_{\bfit m_1,\ldots,\bfit m_{\N-1}}=\prod\limits_{n=1}^{N-1}
\prod\limits_{k=1}^{n}\sigma_{k}^{m_{nk}}(\bgamma_{n})w_\N,
\ee
where $\!\sigma_{k}(\bgamma_{n})\!$ is the elementary symmetric function of
the variables $\!\gamma_{n1},\ldots,\gamma_{nn}\!$ of order $\!k\!$:
\be\label{bas2}
\sigma_k(\bgamma_n)=\sum_{j_1<\ldots<j_k}
\gamma_{nj_1}\ldots\gamma_{nj_k}.
\ee
Similarly, the Whittaker module $W'={\cal U}w_\N'$ is spanned by the
polynomials
\be
w'_{\bfit m_1,\ldots,\bfit m_{\N-1}}=
 \prod\limits_{n=1}^{N-1}
\prod\limits_{k=1}^{n}\sigma_{k}^{m_{nk}}(\bgamma_{n}).
\ee
The Whittaker modules $W$ and $W'$ are irreducible.
\end{te}
{\bf Proof.}
Let us prove the statement for the module $W$. Using the identity
$e^{i\hbar\d_{\gamma_{nj}}}\sigma_k(\bgamma_n)=\sigma_k(\bgamma_n)
+i\hbar\sum\limits_{r=1}^{k-1}(-1)^r
\sigma_{k-r-1}(\bgamma_n)\gamma_{nj}^r$, and the explicit formula for
the action of $E_{n+1,n}$ on the Whittaker vector:
\be\label{low1}
E_{n+1,n}w_\N=i\hbar^{-1}w_\N\sum_{j=1}^n
\frac{\prod\limits_{r=0}^{n+1}
(\gamma_{nj}-\gamma_{n+1,r}+\frac{i\hbar}{2})}
{\prod\limits_{s\neq j}(\gamma_{nj}-\gamma_{ns})}\,,
\ee
one shows, due to (\ref{id2}), that the Whittaker module is
spanned by  elements of the form (\ref{bas1}).

To prove the irreducibility of Whittaker modules we need
the following fact (\cite{Ko2}, Theorem 3.6.2).
{\it Let $M$ be any
${\cal U}$-module which admits an infinitesimal character.
Assume $w\in M$
is a Whittaker vector. Then the submodule ${\cal U}w\subset M$ is
irreducible.}

Proposition \ref{lic} and the above theorem implies that
the modules $W$ and $W'$ are irreducible. \square

Let us note that for any  subalgebra
${\cal U}(\frak{gl}(n))\subset{\cal U}(\frak{gl}(N))$, $2\leq n<N$,
the module over the ring of the polynomials in $\bgamma_{n}$ with the
basis $\prod\limits_{s=1}^{n-1} \prod\limits_{k=1}^{s}
\sigma_{k}^{m_{sk}}(\bgamma_{s})w_\N$ is a ${\cal U}(\frak{gl}(n))$
Whittaker module.

Now one can calculate the explicit form of the action of the central
elements of ${\cal U}(\frak{gl}(n))$ on the space $M$.
It is well known that the generating function $\ca_n(\la)$ of the
central elements  of ${\cal U}(\frak{gl}(n))$ (the Casimir
operators) can be represented as follows
\cite{Weyl}:
\be\label{cas1}
{\cal A}_n(\la) \\
\hspace{-5mm}
=\,\sum_{p\in P_n}{\rm sign}\,p
\Big[(\la-i\hbar\rho^{(n)}_1)\delta_{p(1),1}-
i\hbar E_{p(1),1}\Big]\ldots
\Big[(\la-i\hbar\rho^{(n)}_n)\delta_{p(n),n}-
i\hbar E_{p(n),n}\Big],\hspace{-0.8cm}
\ee
where $\rho^{(n)}_k=\frac{1}{2}(n-2k+1),\;k=1,\ldots, n$ and
the summation is over elements of the permutation group $P_n$.
\begin{prop}
The operators (\ref{cas1}) have the following form  on $M$:
\be\label{cas2}
\hspace{2cm}
\ca_n(\la)=\prod_{j=1}^n(\la-\gamma_{nj}),\hspace{1cm}(n=1,\ldots,N).
\ee
\end{prop}
{\bf Proof.}
By Proposition \ref{lic} it
is sufficient to calculate the action of  $\ca_n(\la)$ on the Whittaker
vector $w_\N$ (or $w_\N'$). Due to (\ref{fww''}), one has
\be\label{char1}
\ca_n(\la)w_\N=\ca_{n-1}(\la-\ts{i\hbar}{2})[\la-i\hbar\rho^{(n)}_n-
i\hbar E_{nn}]w_\N\\
-\,i\hbar\sum_{k=2}^n\ca_{n-k}(\la-\ts{ik\hbar}{2})E_{n,n-k+1}w_\N\,.
\ee
Furthermore, using the relation (\ref{low1}) it is easy to calculate the
action of the composite  generators:
\be\label{low2}
E_{n,n-k+1}w_\N=i\hbar^{-1}w_\N\sum_{j_1=1}^{n-1}\ldots
\sum_{j_k=1}^{n-k+1}\prod_{p=1}^{k-1}
\frac{\prod\limits_{r_p\neq j_{p-1}}(\gamma_{n-p.j_p}-
\gamma_{n-p+1,r_p}+\frac{i\hbar}{2})}
{\prod\limits_{s_p\neq j_p}(\gamma_{n-p,j_p}-
\gamma_{n-p,s_p})}\,.
\ee
By substituting (\ref{low2}) into  (\ref{char1}) one arrives,
after some combinatoric, at (\ref{cas2}). Obviously, we can also do
the analogous calculation for the Whittaker module $W'$.\square

It  remains to construct a pairing between $W$ and $W'$, and to prove that
the Whittaker modules  $W$ and $W'$ are dual with respect to this pairing.
\begin{de} Let $\phi\in W'$ and $\psi\in W$. Define the pairing
$\langle\,,\rangle$: $W'\otimes W\rightarrow\CC$ by
\be\label{in1'}
\langle\phi,\psi\rangle=\int\limits_{\RR^{\frac{N(N-1)}{2}}}
\mu_0(\bgamma)\ov{\phi}(\bgamma)\,\psi(\bgamma)\,
\prod\limits_{\stackrel{\scriptstyle n=1}{j\leq n}}^{N-1}d\gamma_{nj}\,,
\ee
where
\be\label{in2'}
\mu_0(\bgamma)=\prod_{n=2}^{N-1}\prod\limits_{s<p}
(\gamma_{ns}-\gamma_{np})
(e^{\frac{2\pi\gamma_{np}}{\hbar}}-e^{\frac{2\pi\gamma_{ns}}{\hbar}}).
\ee
\end{de}

The integral (\ref{in1'}) converges absolutely. Actually, this is
a simple corollary of a more general statement.  Denote the integrand in
the right hand side of (\ref{in1'}) by $I_\N(\bgamma)$.
\begin{lem}
For $\gamma_{nj}$, such that
$\min\limits_j\{{\rm Im}\,\gamma_{nj}\}>
\max\limits_k\{{\rm Im}\,\gamma_{n+1,k}\}-\frac{\hbar}{2}$,
the following inequality holds:
\be\label{absol}
|I_\N(\bgamma)|\leq
|P(\bgamma)|\,\exp\Big\{-\frac{1}{(2N-3)!!}\sum\limits_{n=1}^{N-1}
\sum\limits_{j=1}^{n}|{\rm Re}\,\gamma_{nj}|\Big\}\,
\ee
for some polynomial $P(\bgamma)$ in $\bgamma_1,\ldots,\bgamma_{\N-1}$.
\end{lem}
{\bf Proof.} For any $x>0$ the inequality
$|\Gamma(x+iy)|\leq\Gamma(x)|p_x(y)|\cosh^{-\frac{1}{2}}(\pi y)$ is valid
for some polynomial $p_x(y)$ with degree depending on
$x$. Hence, the integrand satisfies:
\be
|I_\N(\bgamma)|\leq|P(\bgamma)|\,e^{\frac{\pi}{\hbar}S_N(\bgamma)},
\ee
where
\be
S_\N(\bgamma)=\sum_{n=2}^{N-1}\sum_{s<p}|{\rm Re}\,\gamma_{ns}-
{\rm Re}\,\gamma_{sp}| \\
-\,\frac{1}{2}\sum_{n=1}^{N-2}\sum_{j,k}
|{\rm Re}\,\gamma_{nj}-{\rm Re}\,\gamma_{n+1,k}|-
\frac{N}{2}\sum_{j=1}^{N-1}|{\rm Re}\,\gamma_{\N-1,j}|.
\ee
Using the inequality
\be
\sum_{j,k=1}^m|a_j-b_k|-\sum_{j<k}^m|a_j-a_k|-
\sum_{j<k}^m|b_j-b_k|\geq 0\, ,
\ee
which holds for any set of real parameters
$\{a_j,b_j\};\,j=1,\ldots, m$ (see (11.118)
in \cite{GR}), one proves the recursive relation
$S_\N\leq\frac{1}{2N-3}(S_{\N-1}-\sum_j|{\rm Re}\,\gamma_{\N-1,j}|)\, $.
This implies  (\ref{absol}).\square
\begin{prop}
Let $\bgamma_\N\in\RR^N$. The Whittaker modules $W$ and $W'$ are dual with
respect to the pairing defined by (\ref{in1'}). I.e. for any
$\phi\in W'$ and $\psi\in W,$ the generators $X\in\frak{gl}(N,\RR)$
possess the property
\be\label{inv5}
\langle\phi, X\psi\rangle\,=\,-\,\langle X\phi,\psi\rangle.
\ee
\end{prop}
{\bf Proof.}
For example, consider $\langle\phi,E_{n,n+1}\psi\rangle$, where
$E_{n,n+1}$ is defined by (\ref{m1b}). It is easy to see that the
expression $\ov{E_{n,n+1}\phi}\,\psi$ does not have poles
in the upper hyper-plane and it is therefore possible to deform the
integration contour. Consider the shifts
$\gamma_{nj}\to\gamma_{nj}+i\hbar$ and note the difference equation
\be\label{in3}
e^{i\hbar\d_{\gamma_{nj}}}\mu_0(\bgamma)=\mu_0(\bgamma)\,
\prod_{s\neq j}\frac{\gamma_{nj}-\gamma_{ns}+i\hbar}
{\gamma_{nj}-\gamma_{ns}}\,.
\ee
Then, using a deformation of the contour and the estimates (\ref{absol}),
one obtains (\ref{inv5}).
\square

\section{The orbit method for $GL(N,\RR)$ and the Gelfand-Zetlin
representation}

In this section we discuss the orbit method  interpretation
of  the explicit algebraic construction of the
infinite-dimensional representation
of ${\cal U}(\frak{gl}(N))$ given in the previous  section.
Similar results for the case
of a compact group were obtained  in \cite{AFS}, \cite{S}.

Let $G$ be a Lie group, $\frak{g}$ the corresponding Lie algebra, and let
 $\frak{g}^*$ be its dual. The
coadjoint orbits $\mathcal{O}$ of $G$ are equipped
with a canonical symplectic structure  given by the Kirillov-Kostant
two-form \cite{Kir}. The space of linear functions
on $\frak{g}^*$  is closed under the Kirillov-Kostant Poisson bracket and
coincides with $\frak{g}$ as a Lie algebra.
Consider the  coadjoint orbits of
$GL(N,\RR)$ with  stabilizer $H_\N$ $=GL(1,\RR)^{\otimes N}$.
Below,  we construct Darboux coordinates on the finite cover
$\widetilde{\mathcal{O}}^{(0)}$ of the open part
$\mathcal{O}^{(0)}$ of such  an orbit  $\mathcal{O}$.
We also give expressions in these coordinates for the restriction
onto  $\mathcal{O}^{(0)}$ of the linear
functions on  $\frak{g}^*$. The  corresponding Hamiltonian vector fields
act on $\mathcal{O}^{(0)}$. The explicit realization
of ${\cal U}(\frak{gl}(N))$
by difference operators given in the previous  section
may be considered as a ``quantization'' of this Poisson algebra.
Let us stress that this construction gives a representation
of ${\cal U}(\frak{gl}(N))$ which cannot be integrated to a
representation of the  group. Our discussion below  shows that this
is quite natural because  part
of the orbit
$\mathcal{O}^{(0)}$ is not stable with respect to the action
of the group on the  orbit $\mathcal{O}$.


Let $G=GL(N,\RR)$ and $\frak{g}=\frak{gl}(N,\RR)$. We identify
$\frak{g}$ and its dual $\frak{g}^*$ via the Killing
form. Consider the  coadjoint orbit $\mathcal{O}$
of the diagonal element $u_\N\in \frak{g}^*$ with
pairwise different diagonal entries:
$\bgamma_\N=(\g_{\N 1},\g_{\N 2},\ldots,\g_{\N\N})$:
\bqa
u_\N=
\begin{pmatrix}
\g_{\N 1} & \ldots & 0  \\
\vdots & \ddots  &\vdots \\
0 & \cdots & \g_{\N\N}
\end{pmatrix}.
\eqa
The stabilizer of $u_\N$ is $H_{\N}=GL(1,\RR)^{\otimes \N}$and the
general element of the orbit may be parameterized as
$u=g^{-1}u_\N g$, where $g \in GL(N,\RR)$ is defined up to the
left action of $H_\N$. There is a canonical
symplectic two-form on $\mathcal{O}$, given
in this parameterization by
\be
\Omega=Tr\, u_\N (dg\, g^{-1})^2\,.
\ee
Consider the subspace in $GL(N,\RR)$ which
consists of the elements:
\bqa
\label{ge}
g=g_\N g_{\N-1}\cdots g_2\,.
\label{mult}
\eqa
The matrices $g_n$ differ from the unit matrix only in the upper-left
$n\times n$ corner and are defined recursively as follows.
Denote the upper-left $(n\times n)$ sub-matrix of $g_n$ by $f_n$.
There should exist certain $\g_{n-1,j}\,,\,j=1,\ldots,n-1$ so that
the following condition holds:
\bqa \label{first}
f^{-1}_n
\begin{pmatrix}
\g_{n1} & 0 & \cdots & 0 \\
0 & \ddots &  & \vdots \\
\vdots &  & \g_{n,n-1} & 0 \\
0 & \cdots & 0 & \g_{nn}
\end{pmatrix}
f_n=\begin{pmatrix}
\g_{n-1,1} & 0 & \cdots & * \\
\vdots  & \ddots & 0 & \vdots \\
0 & \cdots  & \g_{n-1,n-1} &* \\
* & \cdots & * & *
\end{pmatrix}.
\eqa
This means that the diagonal matrix conjugated by $f_n$ can have non-zero,
non-diagonal entries only in the last column and on the last row.
The matrix $f_n$ is defined uniquely by $\g_{nj},\g_{n-1,k}$ up to
the left action by the diagonal matrix and the right action by
the diagonal matrix with the last diagonal element being equal to
one. Let us choose the following representatives for $f_n$ in the coset
of $G$ by the left action by  diagonal matrices:
\be\label{repr}
\hspace{1.8cm}
(f_n)_{jk}=\frac{Q_{n-1,k}}{\gamma_{nj}-\gamma_{n-1,k}}
\frac{\prod\limits_{r=1}^{n-1}
(\gamma_{nj}-\gamma_{n-1,r})}
{\prod\limits_{s\neq j}(\gamma_{nj}-\gamma_{ns})}
\hspace{1cm}(k=1,\ldots,n-1),\\
\hspace{-3cm}
(f_n)_{jn}=\frac{\prod\limits_{r=1}^{n-1}(\gamma_{nj}-\gamma_{n-1,r})}
{\prod\limits_{s\neq j}(\gamma_{nj}-\gamma_{ns})}\;,
\hspace{-3cm}
\ee
where $Q_{n-1,k}$, $k=1,\ldots,n-1$ are additional coordinates on the orbit.
This gives a  map from $\RR^{N(N-1)}\backslash\Delta$
to the open part $\mathcal{O}^{(0)}$ of the orbit. Here,
$\Delta$ is a union of three subspaces
of $\RR^{N(N-1)}$: $\Delta=D_1\cup D_2\cup D_3$.
The subspace $D_1$ is a subspace where at least
two $\bgamma$-coordinates of the same level coincide: $\g_{nj}=\g_{nk}$.
The subspace $D_2$ is defined  as a subspace where
at least two of the $\bgamma$-coordinates
from  consecutive levels coincide
$\g_{nj}=\g_{n-1,k}$. Finally, the subspace $D_3$ is  a subspace
where at least one  of $Q_{nj}$ is  zero.
Now,   let us construct  the inverse
map. Consider the element of the
orbit $u=g^{-1}u_\N g$ shifted by the $N\times N$ unit matrix
$I$ multiplied by a formal variable $\la$:
\bqa
T(\lambda)=\la I-u\,.
\eqa
The $n$-th principal
minor $a_n(\la)$ of $T(\lambda)$ does not depend on $g_m$,  $m<n+1$.
For the general point of the orbit, the roots of the minors may be
complex. Consider the open  part $\mathcal{O}^{(0)}$
of the orbit $\mathcal{O}$  such that all upper-left sub-matrices
are diagonalizable and their determinants have distinct
real roots. We also impose the
condition that there are no identical  roots of  minors
with consecutive ranks. Then,  define the $\bgamma$-variables
to be the roots of the minors:
\bqa\label{a1}
a_n(\lambda)=\prod_{j=1}^n(\lambda-\g_{nj})\,.
\eqa
Note that, by this condition
the roots are defined only up to the action of the symmetric
group $S_n$. As a result we obtain
the coordinates on the finite cover of an open part of the orbit.

Now, we give the explicit expressions for the coordinates $Q_{nj}$.
Let $b_n(\la)$
be the determinant of the $n\times n$ corner sub-matrix of the matrix
 obtained from $T(\lambda)$ by interchanging the
$n$-th and the $(n+1)$-th columns. I.e.
$b_n(\la)=\det\limits_{n\times n}||T(\la)S_{n,n+1}||$, where the
$N\times N$ matrix $S_{n,n+1}$ is a unit matrix
with the diagonal
$2\times 2$ sub-matrix in the $n$, $n+1$ columns and the
$n$, $n+1$ rows  replaced by the matrix
$\left(\begin{array}{cc} 0 & 1\\ 1 & 0
\end{array}\right)$.
Similarly, we define $c_n(\la)$ as the  determinant
of the $n\times n$
corner sub-matrix obtained from $T(\la)$ by interchanging the $n$-th
and the $(n+1)$-th rows; i.e.
$c_n(\la)=\det\limits_{n\times n}||S_{n,n+1}T(\la)||$.
Then, the rest of the variables may be introduced as
follows:
\be Q_{nj}=b_n^{-1}(\g_{nj})a_{n+1}(\g_{nj})=-
c_n(\g_{nj}) a^{-1}_{n-1}(\g_{nj})\,,
\ee
where  the last equality can be proved by a straightforward
check.
It is not difficult to verify that this definition is compatible with
the explicit parameterization in terms of $(\g_{nj},Q_{nj})$
given above. Due
to the condition (\ref{first}), the $n$-th left-upper corner
sub-matrix of the matrix $u_n=g^{-1}_{n+1}\cdots g^{-1}_{\N}u_\N
g_\N \cdots g_{n+1}$  is  diagonal matrix;
$a_n(\lambda)$ is given by (\ref{a1}).
The determinant
$b_n(\la)$ is invariant under conjugation by the matrix which
differs from the unit matrix only in the upper-left
$(n-1)\times(n-1)$ corner. Thus, $b_n(\la)$ may be represented
as:
\be
b_n(\lambda)=\det(g_n^{-1})\,
\det_{n\times n}||(\la-u_n)g_nS_{n,n+1}||\,.
\ee
A straightforward calculation gives
\be\label{b1}
b_n(\la)=\sum_{j=1}^n Q^{-1}_{nj}
\prod\limits_{r=1}^{n+1}(\gamma_{nj}-\gamma_{n+1,r})
\prod_{s\neq j}\frac{\la-\gamma_{ns}}{\gamma_{nj}-\gamma_{ns}}\,.
\ee
Similarly, for  $c_n(\la)$ one has
\be\label{c1}
c_n(\la)=-\sum_{j=1}^n Q_{nj}
\prod\limits_{r=1}^{n-1}(\gamma_{nj}-\gamma_{n-1,r})
\prod_{s\neq j}\frac{\la-\gamma_{ns}}{\gamma_{nj}-\gamma_{ns}}\,.
\ee
This gives the inverse transformation
from the finite cover of the open part of the orbit
$\widetilde{\mathcal{O}}^{(0)}$ to $\RR^{\frac{N(N-1)}{2}}\backslash\Delta$.

\hspace{-0.3cm}
A simple calculation of the Kirillov-Kostant symplectic form
in the coor\-di\-na\-tes $\!(\g_{nj},Q_{nj})\!$
gives the following:
\begin{prop}
\begin{enumerate}
\item The coordinates $(\g_{nj},Q_{nj})$  parameterize the finite
cover $\widetilde{\mathcal{O}}^{(0)}$ of
the open part of the orbit $\mathcal{O}^{(0)}$. The product
of the symmetric groups $S=\prod_{n=2}^{\N}  S_n$ acts freely on the
fibres of the projection: $\pi:
\widetilde{\mathcal{O}}^{(0)}\rightarrow \mathcal{O}^{(0)}$
and the factor is $\mathcal{O}^{(0)}$.

\item The lift of the 2-form $\Omega$ on
 $\widetilde{\mathcal{O}}^{(0)}$ has the canonical form
\be\label{al}
\pi^* \Omega=\sum_{n=1}^{N-1}\sum_{j=1}^nd
\g_{nj}\wedge Q_{nj}^{-1}dQ_{nj} .
\ee
The coordinates  $(\g_{nj}, Q_{nj})$ give  Darboux
coordinates on the finite cover
of the open part $\mathcal{O}^{(0)}$ of the coadjoint orbit of
$GL(N,\RR)$ with the stabilizer conjugated to
$H_\N=GL(1,\RR)^{\otimes N}$.
\end{enumerate}
\end{prop}


Now let us give explicit  expressions for the  elements
on the super-diagonal, diagonal and over-diagonal
of the  matrix
$u=g^{-1}u_\N g$ in the GZ-coordinates.
These give the classical counterparts of the generators (\ref{m1}).
\begin{prop}
The  linear functions $u_{nn},\,u_{n,n+1}\,,
u_{n+1,n}$  on
$\frak{gl}(N,\RR)^*$, being
restricted to $\mathcal{O}^{(0)}$ and lifted by the projection map
$\pi:\widetilde{\mathcal{O}}^{(0)}\rightarrow \mathcal{O}^{(0)}$
to $\widetilde{\mathcal{O}}^{(0)}$, have the following
form in  GZ-coordinates
\footnote{There is another representation with
``square roots'' which is  closer to that  given in \cite{GelZ}.
It may be obtained by a different choice of  parameterization for
the solution of (\ref{first}).
}
:
\be\label{gez1}
u_{nn}=\sum_{j=1}^n\g_{nj}-\sum_{j=1}^{n-1}\g_{n-1,j}\,,\\
 u_{n,n+1}=-\sum_{j=1}^n
\frac{\prod\limits_{r=1}^{n+1}(\gamma_{nj}-\gamma_{n+1,r})}
{\prod\limits_{s\neq j}(\gamma_{nj}-\gamma_{ns})}\;
Q_{nj}^{-1}\,,\\
u_{n+1,n}=\sum_{j=1}^n
\frac{\prod\limits_{r=1}^{n-1}(\gamma_{nj}-\gamma_{n-1,r})}
{\prod\limits_{s\neq j}(\gamma_{nj}-\gamma_{ns})}\;Q_{nj}\,.
\ee
\end{prop}
{\bf Proof.} Let us start with the diagonal elements. The coadjoint
action of the
diagonal matrix $\Phi={\rm diag}(e^{\phi_1}, \cdots e^{\phi_\N})$
on $u$
may be represented as the right action on the element
$g$, (\ref{ge}). It is easy to see that this is equivalent to
a shift of the variables:
\be
Q_{nj} \rightarrow Q_{nj}e^{\phi_n-\phi_{n+1}},\\
Q_{\N j} \rightarrow  Q_{\N j}e^{\phi_\N}.
\ee
Thus, the corresponding Hamiltonian functions $u_{nn}$ are:
\be
u_{nn}=\sum_{j=1}^n\g_{nj}-\sum_{j=1}^{n-1}\g_{n-1,j}\,.
\ee
This result could also be derived from the explicit
representation of the matrix $u$.
Moreover, from the structure of the matrices
$u_n, f_n$ one could easily obtain expressions for
the other generators. For example,
\be\label{low}
u_{n+1,n}=\sum_{j=1}^n(u_n)_{n+1,j}(f_{n})_{jn}=
\sum_{j=1}^{n}
\frac{\prod\limits_{r=1}^{n-1}(\gamma_{nj}-\gamma_{n-1,r})}
{\prod\limits_{s\neq j}(\gamma_{nj}-\gamma_{ns})}\,
Q_{nj}\,.
\ee
Quite similarly,
\be\label{raise}
u_{n,n+1}=-\sum_{j=1}^{n}
\frac{\prod\limits_{r=1}^{n+1}(\gamma_{nj}-\gamma_{n+1,r})}
{\prod\limits_{s\neq j}(\gamma_{nj}-\gamma_{ns})}\,
Q_{nj}^{-1}\,.
\ee
\square

There is  a representation for
$u_{n,n},\,u_{n,n+1},$ and $u_{n+1,n}$ in terms of contour integrals
using the polynomials
$a_n(\lambda),b_n(\lambda),$ and $c_n(\lambda)$:
\be \label{integr}
u_{n,n}=-\frac{1}{2\pi i}\oint\frac{{ a}_n(\la)}
{{ a}_{n-1}(\la)}\frac{d\la}{\la}\,, \\
u_{n,n+1}=
-\frac{1}{2\pi i}\oint{ a}_n^{-1}(\la){ b}_n(\la)d\la\,,\\
u_{n+1,n}=
-\frac{1}{2\pi i}\oint{ c}_n(\la){ a}_n^{-1}(\la)d\la\,.
\ee
In the next section we discuss the connection of
this representation with  Drinfeld's  ``new realization''
of the Yangian \cite{Dr2}.

\section{The  Yangian $\!\!Y(\frak{gl}(N))\!$ and the Gelfand-Zetlin
representation}

In the previous section we outlined the construction of the
Gelfand-Zetlin parameterization of the
coadjoint orbit of $GL(N,\RR)$. The
polynomials $a_n(\la),b_n(\la)$ and $c_n(\la)$, which appear in the
invariant formulation of this parameterization, bear an obvious
similarity to the basic ingredients of Drinfeld's ``new realization''
of the Yangian. In this section, we derive the explicit expressions
for the Drinfeld generators of the Yangian $Y(\frak{gl}(N))$
in the Gelfand-Zetlin realization\footnote{For the
 connection of the finite-dimensional Gelfand-Zetlin
representation
with the Yangian see \cite{Cher1}, \cite{NT1}.
}. These explicit expressions will be important for
a new interpretation
of the QISM approach  to the solution of the open Toda chain
\cite{KL2}.

We recall some well-known facts about the Yangian
$Y(\frak{gl}(N))$, \cite{Dr2} (see also the recent review \cite{MNO}).
The Yangian $Y(\frak{gl}(N))$ is an associative Hopf
algebra generated by the elements $T_{ij}^{(r)}$, where
$i,j=1,\ldots,N$ and $r=0,\ldots,\infty$, subject the following
relations. Consider the $N\times N$ matrix
$T(\la)=||T_{ij}(\la)||_{i,j=1}^\N$ with  operator valued
entries
\be\label{y1}
T_{ij}(\la)=\la\delta_{ij}+\sum_{r=0}^\infty T_{ij}^{(r)}\la^{-r}.
\ee
Let
\be\label{r-mat}
\hspace{2cm}
R_\N(\la)\,=\,I\otimes I+i\hbar P/\la\;;\hspace{1cm}
P_{ik,jl}=\delta_{il}\delta_{kj}
\ee
be an $N^2\times N^2$ numerical matrix (the Yang $R$-matrix).
Then the relations between the generators $T_{ij}^{(r)}$ can be
written in the standard form
\be\label{y2}
R_\N(\la-\mu)(T(\la)\otimes I)(I\otimes T(\mu))=
(I\otimes T(\mu))(T(\la)\otimes I)R_\N(\la-\mu).
\ee
The centre of the Yangian is generated by the coefficients of the
following formal Laurent series
(the quantum determinant of $T(\la)$ in the sense of \cite{KS}):
\be\label{y4}
{\rm det}_q T(\la) \\ =\,\sum_{p\in P_\N}{\rm sign}\,p\;
T_{p(1),1}(\la-i\hbar\rho^{(\N)}_1)\ldots
T_{p(k),k}(\la-i\hbar\rho^{(\N)}_k)
\ldots T_{p(\N),\N}(\la-i\hbar\rho^{(\N)}_\N)\,,
\hspace{-0.5cm}
\ee
where $\rho^{(\N)}_n=\ts{1}{2}(N-2n+1),\,(n=1,\ldots,N)$ and the
summation is  over elements of the permutation group
$P_\N$.
Let $X(\la)=||X_{ij}(\la)||_{i,j=1}^n$ be an $n\times n$
sub-matrix of the matrix $||T_{ij}(\la)||_{i,j=1}^\N$.
It is obvious from the explicit form of $R_\N(\la)$
that this sub-matrix satisfies an analogue of the relations
(\ref{y2}).
The quantum determinant ${\rm det}_q X(\la)$ is defined similarly
to (\ref{y4}) (with the  evident change $N\to n$).

The following way to describe the Yangian $Y(\frak{gl}(N))$
was introduced in \cite{Dr2}.
Let ${\bf A}_n(\la)$, $n=1,\ldots, N$, be the quantum
determinants of the sub-matrices, determined by  the first $n$
rows and columns, and let  the
operators ${\bf B}_n(\la),{\bf  C}_n(\la)$,
$n=1,\ldots, N-1$,  be the quantum
determinants of the  sub-matrices with  elements $T_{ij}(\la)$,
where $i=1,\ldots,n$; $j=1,\ldots,n-1,n+1$ and
$i=1,\ldots,n-1,n+1$; $j=1,\ldots,n$, respectively.
The expansion coefficients of
$ {\bf A}_n(\la),{\bf B}_n(\la) ,{\bf C}_n(\la)$,
$n=1,\ldots,N-1$, with respect to $\la$,
together with those of ${\bf A}_\N(\la)$, generate the algebra
$ Y(\frak{gl}(N)) $.
The parts of the relations between the Drinfeld  generators,
which we use below, are as follows:
\be\label{cw1}
\hspace{3cm}
[{\bf A}_n(\la),{\bf A}_m(\mu)]=0\;;\hspace{1cm}(n,m=1,\ldots, N),\\

\hspace{2cm}[{\bf B}_n(\la),{\bf B}_m(\mu)]=0\,;\;\ \
[{\bf C}_n(\la),{\bf C}_m(\mu)]=0\;;\hspace{1cm} (m\neq n\pm 1),\\
(\la-\mu+i\hbar){\bf A}_n(\la){\bf B}_n(\mu)\,=\,
(\la-\mu){\bf B}_n(\mu){\bf A}_n(\la)+i\hbar{\bf A}_n(\mu){\bf B}_n(\la),\\
(\la-\mu+i\hbar){\bf A}_n(\mu){\bf C}_n(\la)\,=\,
(\la-\mu){\bf C}_n(\la){\bf A}_n(\mu)+i\hbar{\bf A}_n(\la){\bf C}_n(\mu).
\ee

Let $A(\frak{gl}(N))$ be the commutative subalgebra of
$Y(\frak{gl}(N))$ generated by ${\bf A}_n(\la),\,n=1,\ldots,N$.
It was proved in \cite{Cher2} that $A(\frak{gl}(N))$ is
the maximal commutative subalgebra of $Y(\frak{gl}(N))$.

There is a natural  epimorphism
$\pi_\N :Y(\frak{gl}(N))\rightarrow{\cal U}(\frak{gl}(N))$,
\be\label{bz1}
\hspace{3cm}
\pi_{\N}(T_{jk}(\la))=\la\delta_{jk}-i\hbar E_{jk}\,,
\hspace{1cm}(j,k=1,\ldots, N).
\ee
Denote the images under $\pi_{\N}$  of the
 generators ${\bf A}_n(\la) ,{\bf B}_n(\la)$,
and ${\bf C}_n(\la)$
 by
$\ca_n(\la),\cb_n(\la)$ and $\cc_n(\la)$,   respectively.
Let us describe the images of the Drinfeld  generators.
The images of the generators of $A(\frak{gl}(N))$
under the homomorphism (\ref{bz1}) have the form (\ref{cas1}).
In particular, the image of ${\bf A}_{\N}(\la)$, which is a polynomial
in ${\la}$ of order $N\!$, gives a generating function for the
generators of the centre ${\cal Z}\subset{\cal U}(\frak{gl}(N))$.
For the other generators $\cb_n(\la)$, $\cc_n(\la)$ we have
\begin{lem}
\be\label{bzb}
\cb_n(\la)=[\ca_n(\la),E_{n,n+1}],
\ee
\be\label{bzc}
\cc_n(\la)=[E_{n+1,n},\ca_n(\la)],
\ee
where $n=1,\ldots,N-1$.
\end{lem}\
{\bf Proof.} A direct computation using (\ref{y4}) and the explicit
expressions for the images of the generators
${\bf A}_n(\la) ,{\bf B}_n(\la)$
and ${\bf C}_n(\la)$ under the homomorphism $\pi_{\N}$.
\square

Using (\ref{bzb}), (\ref{bzc}) and the explicit expressions
for the generators $E_{jk}$ introduced in Proposition 2.1
we arrive at  the following
\begin{prop}\label{boz}
The operators
\be\label{bza}
\ca_{n}(\lambda)=\prod\limits_{j=1}^{n}(\la-\gamma_{nj}),\\
 \cb_{n}(\la)=\sum_{j=1}^{n}\prod_{s\neq j}
\frac{\la-\gamma_{ns}}{\gamma_{nj}-\gamma_{ns}} \,
\prod_{r=1}^{n+1}(\gamma_{nj}-\gamma_{n+1,r}-\frac{i\hbar}{2})
\,e^{-i\hbar\d_{\gamma_{nj}}}, \\
\cc_{n}(\la)=-\sum_{j=1}^{n}\prod_{s\neq j}
\frac{\la-\gamma_{ns}}{\gamma_{nj}-\gamma_{ns}} \,
\prod_{r=1}^{n-1}(\gamma_{nj}-\gamma_{n-1,r}+\frac{i\hbar}{2})
\,e^{i\hbar\d_{\gamma_{nj}}}
\ee
define a  representation of the Yangian.
\end{prop}

We notice that the explicit expressions
for $\ca_n(\la) ,\cb_n(\la),\cc_n(\la)$
may be obtained directly from the  defining relations of the Yangian.
One may start with the following
representation of the maximal commutative subalgebra
$A(\frak{gl}(N))$ by polynomials with real roots:
$\ca_{n}(\lambda)=\prod\limits_{j=1}^{n}(\la-\gamma_{nj}),\,
n=1\ldots,N$. Then, we should resolve the rest of the  Yangian
relations to find the explicit expressions for the generators
$\cb_n(\la)$ and $\cc_n(\la)$
in terms of some operators acting on the space of functions depending
on the variables $\gamma_{nj}\,,j=1,\ldots, n;\,n=1,\ldots,N-1$.
The full set of Yangian relations in terms of ${\bf A}_n(\la),
{\bf B}_n(\la),{\bf  C}_n(\la)$
is known implicitly through the ``new realization''
by Drinfeld \cite{Dr2} (example after Theorem 1).
This may be used to fix
${\bf B}_n(\la), {\bf C}_n(\la)$ up to conjugation
by an arbitrary function of $\gamma_{nj}$. The  expressions
in (\ref{bza}) are obtained by choosing the appropriate representative.

The following integral formulas (which  may be considered as a
``quantization'' of the classical relations (\ref{integr}))
express the generators of ${\cal U}(\frak{gl}(N))$ in terms of
$\ca_n(\la),\cb_n(\la),\cc_n(\la)$:
\be\label{rtt3}
E_{n,n+1}=\frac{1}{2\pi \hbar}\oint
{\cal A}_n^{-1}(\la){\cal B}_n(\la) d\la\,,\\
E_{n+1,n}=\frac{1}{2\pi\hbar}\oint {\cal C}_n(\la)
{\cal A}_n^{-1}(\la) d\la\,,\\
E_{nn}=\frac{1}{2\pi\hbar}
\oint\limits \frac{{\cal A}_n(\la)}{{\cal A}_{n-1}
(\la-\ts{i\hbar}{2})}\frac{d\la}{\la}-\frac{1}{2}(n-1)\,.
\ee
Here, the integrands are understood as Laurent series and the
contours of integrations are taken around $\infty$.
This representation is similar to the expressions of the generators
in the ``new realization'' of the Yangian through
${\bf A}_n(\la),{\bf B}_n(\la),{\bf  C}_n(\la)$ \cite{Dr2}
(example after Theorem 1).


\section{Application to  quantum integrable systems}
\subsection{The open Toda chain}
The open Toda chain  corresponding to $\frak{gl}(N,\RR))$
is one of the simplest
examples of an integrable quantum mechanical system with $N$ degrees
of freedom (see \cite{STS} and references therein).
It has $N$ mutually commuting Hamiltonians, the first two of which
are given by
\be\label{ham}
h_1=\sum\limits_{j=1}^N p_j\,,\\
h_2=\sum_{j<k}p_jp_k-\sum_{j=1}^{N-1}e^{x_{j}-x_{j+1}}\,,
\ee
where $[x_n,p_m]=i\hbar\delta_{nm}$.

Recall the general idea behind the Representation Theory
approach to  quantum integrability in the case of the open Toda
chain. Let $V$ and $V'$ be
any dual irreducible Whittaker ${\cal U}(\frak{gl}(N,\RR))$- modules
and $w \in V$ and $w' \in V'$ be the corresponding
cyclic Whittaker vectors with
respective characters $\!\chi_+(E_{n,n+1})=-i\hbar^{-1}\!$
and $\!\chi_-(E_{n+1,n})=-i\hbar^{-1}\!$.
We assume that
the  action of the Cartan subalgebra is integrated to the action of
the Cartan torus,  so that the following function  is well defined
\be\label{pair}
\psi_{\gamma_{\N1},\ldots,\gamma_{\N\N}}=
e^{-\bfit x\cdot\brho^{(N)}}
\langle  w'_\N,e^{-\sum\limits_{k=1}^Nx_kE_{kk}}w_\N \rangle\,.
\ee
This is nothing but the $GL(N,\RR)$ Whittaker function \cite{Jac}
written in terms of the Gauss decomposition \cite{GKM}.
In the case of irreducible Whittaker modules,
the  action of the elements of the centre ${\cal Z}$ of
the universal enveloping algebra ${\cal U}(\frak{gl}(N,\RR))$
is proportional to the action of the unit operator \cite{Ko2} (Theorem
3.6.1).
Consider the function
\be\label{pair1}
\widetilde{\psi}_{\gamma_{\N1},\ldots,\gamma_{\N\N}}=
e^{-\bfit x\cdot\brho^{(N)}}
\langle  w'_\N,e^{-\sum\limits_{k=1}^Nx_kE_{kk}}z\,  w_{\N} \rangle\,,
\ee
where $z$ belongs to the centre ${\cal Z}$.
Using the properties of the Whittaker vectors one could show
that there is a differential operator
$\mathcal{D}_z$ in the variables $x_k$ such that
$\widetilde{\psi}_{\gamma_{\N1},\ldots,\gamma_{\N\N}}=
\mathcal{D}_z\psi_{\gamma_{N1},\ldots,\gamma_{NN}}.$
Thus, taking the first two
Casimir operators, one gets the first two Hamiltonians (\ref{ham}).
On the other hand, $z$ is an element of the centre ${\cal Z}$ and
$\widetilde{\psi}_{\gamma_{N1},\ldots,\gamma_{NN}}$ is proportional
to $\psi_{\gamma_{N1},\ldots,\gamma_{NN}}$ with some
numerical coefficient. Hence, we get the common eigenfunction for the set
of differential operators corresponding to the elements of ${\cal Z}$.

Let us give  an explicit realization of the representation
of ${\cal U}(\frak{gl}(N,\RR))$ in terms of some difference/differential
operators, which integrates to the representation of the Cartan group.
This will lead to the integral formula for the wave function of the open
Toda chain. To find such a formula for the GZ-representation, we substitute
the expressions (\ref{wv}), (\ref{wv'}), and (\ref{m1a}) into
(\ref{pair}), thus obtaining
\be\label{wf1}
\psi_{\gamma_{\N1},\ldots,\gamma_{\N\N}}\;=
e^{-\bfit x\cdot\brho^{(N)}} \\
\hspace{-1cm}
\times\!\int\limits_{\RR^{\frac{N(N-1)}{2}}}
\prod_{n=1}^{N-1}
\frac{\prod\limits_{k=1}^n\prod\limits_{m=1}^{n+1}
\hbar^{\frac{\gamma_{nk}-\gamma_{n+1,m}}{i\hbar}+\frac{1}{2}}\;
\Gamma(\frac{\gamma_{nk}-\gamma_{n+1,m}}{i\hbar}+\frac{1}{2})}
{\prod\limits_{s<p}\left|
\Gamma(\frac{\gamma_{ns}-\gamma_{np}}{i\hbar})\right|^2}
\;e^{\frac{i}{\hbar}
\sum\limits_{n,j=1}^N(\gamma_{nj}-\gamma_{n-1,j})x_n}
\prod_{\stackrel{\scriptstyle n=1}{j\leq n}}^{N-1}d\gamma_{nj}
\,,\hspace{-1.5cm}
\ee
where by definition $\gamma_{nj}=0$ for $j>n$.

In the study of the analytic properties
of this solution with respect to $\bgamma_\N$, the following reformulation
is very useful:
\begin{te}
 An  analytical continuation of the eigenfunction (\ref{wf1}) as a
function of $\bgamma_\N$ can be expressed in the form
\be\label{wf6}
\psi_{\bgamma_{\N}}(x_1,\ldots,x_\N)\\
=\;\int\limits_{\cal S}\prod_{n=1}^{N-1}
\frac{\prod\limits_{k=1}^n\prod\limits_{m=1}^{n+1}
\hbar^{\frac{\gamma_{nk}-\gamma_{n+1,m}}{i\hbar}}\;
\Gamma(\frac{\gamma_{nk}-\gamma_{n+1,m}}{i\hbar})}
{\prod\limits_{s\neq p}
\Gamma(\frac{\gamma_{ns}-\gamma_{np}}{i\hbar})}
\;e^{\frac{i}{\hbar}
\sum\limits_{n,j=1}^N(\gamma_{nj}-\gamma_{n-1,j})x_n}
\prod_{\stackrel{\scriptstyle n=1}{j\leq n}}^{N-1}d\gamma_{nj}\,,
\ee
where the domain of integration  $\!{\cal S}\!$ is defined by the conditions
$\!\min\limits_{j}\{{\rm Im}\,\gamma_{kj}\}>
\max\limits_m\{{\rm Im}\,\gamma_{k+1,m}\}\!$
for all $ k=1,\ldots,N-1 $.
The integral (\ref{wf6}) converges absolutely.
\end{te}
{\bf Proof. }
Let us change the variables of integration in (\ref{wf1}):
\be\label{wf2}
\hspace{3cm}
\gamma_{nj}\;\to\;\gamma_{nj}-\frac{i\hbar}{n}
\sum_{s=1}^n\rho^{(\N)}_s,\hspace{1cm}(n=1,\ldots,N-1)\,,
\ee
and recall that $\rho^{(\N)}_s=\frac{1}{2}(N-2s+1)$.
After elementary calculations, the integral in (\ref{wf1})
acquires the form (\ref{wf6}). It is worth mentioning  that after
the change of variables (\ref{wf2}) we shift  the domain of integration
$\RR^{\frac{N(N-1)}{2}}$ to the complex plane in such a way that
the domain of integration over the variables $\gamma_{n-1,j}$ lies
above the domain of integration over the variables $\gamma_{nj}$.
Thus, we arrive at the analytic  continuation
of the wave function described in the theorem. \square
\begin{ex} Let $N=2$. In this case the
general formula (\ref{wf6}) acquires the form
\be\label{e2}
\psi_{\gamma_{21},\gamma_{22}}(x_1,x_2)\,=\,
\hbar^{-\frac{\gamma_{21}+\gamma_{22}}{i\hbar}}\,
e^{\frac{i}{\hbar}(\gamma_{21}+\gamma_{22})x_2}\\
\times\;
\int\limits_{i\sigma-\infty}^{i\sigma+\infty}
\hbar^{\frac{2\gamma_{11}}{i\hbar}}\,
\Gamma\Big(\frac{\gamma_{11}\!-\!\gamma_{21}}{i\hbar}\Big)
\,\Gamma\Big(\frac{\gamma_{11}\!-\!\gamma_{22}}{i\hbar}\Big)
e^{\frac{i}{\hbar}\gamma_{11}(x_1-x_2)}\,d\gamma_{11}\,,
\ee
where $\sigma\,>\,{\rm max}\{{\rm Im}\,\gamma_{21},
{\rm Im}\,\gamma_{22}\}$. Using the formulas 7.2(15) and
7.12(34) from \cite{BE}, one obtains the solution for the
$GL(2,\RR)$ Whittaker function in terms of the Macdonald function
$K_\nu(z)$:
\be\label{lm1}
\psi_{\gamma_{21},\gamma_{22}}(x_1,x_2)=4\pi\hbar\,
e^{\frac{i}{2\hbar}(\gamma_{21}+\gamma_{22})(x_1+x_2)}
K_{\frac{\gamma_{21}-\gamma_{22}}{\scriptstyle i\hbar}}
\Big(\frac{2}{\hbar}\,e^{(x_1-x_2)/2}\Big)\,.
\ee
\end{ex}

The explicit form of the eigenfunction and the absolute convergence
of the integral lead to a recursive relation
between  the wave functions corresponding to $\frak{gl}(n)$ and
$\frak{gl}(n-1)$; this relation is  a direct consequence of the
Gelfand-Zetlin inductive procedure:
\begin{cor}
Let us fix the solution $\psi_{\gamma_{11}}(x_1)=e^{\frac{i}
{\hbar}\gamma_{11}x_{1}}$ of the $GL(1,\RR)$ Toda chain.
For any $n=2,\ldots,N$  there is a  recursive  relation between
the wave functions $\psi_{\bgamma_n}$ and $\psi_{\bgamma_{n-1}}$:
\be\label{rec}
\psi_{\bgamma_{n}}(x_1,\ldots ,x_n) \\
\hspace{-1.6cm}
=\int\limits_{{\cal S}_n}
\mu^{(n-1)}(\bgamma_{n-1})f_{n}(\bgamma_{n},\bgamma_{n-1})\,
e^{\frac{i}{\hbar}(\sum\limits_{j=1}^{n}
\gamma_{nj}-\sum\limits_{k=1}^{n-1}\gamma_{n-1,k})x_{n}}
\psi_{\bgamma_{n-1}}
(x_1,\ldots ,x_{n-1})
\prod\limits_{ j=1}^{n-1}d\gamma_{n-1,j}\,,
\hspace{-2cm}
\ee
where
\be
\mu^{(n-1)}(\bgamma_{n-1})=\prod_{s\neq p}
\Big[\Gamma\Big(\frac{\gamma_{n-1,s}-\gamma_{n-1,p}}{i\hbar}\Big)
\Big]^{-1},
\ee
and
\be\label{fourc}
f_{n}(\bgamma_{n},\bgamma_{n-1})=
\prod_{k=1}^{n-1}\prod_{m=1}^n
\hbar^{\frac{\gamma_{n-1,k}-\gamma_{nm}}{i\hbar}}
\Gamma\Big(\frac{\gamma_{n-1,k}-\gamma_{nm}}{i\hbar}\Big).
\ee
The domain  of integration ${\cal S}_n$ is defined
by the conditions  $\min\limits_j\{{\rm Im}\gamma_{n-1,j}\}\!>\!
\max\limits_m\{{\rm Im}\gamma_{nm}\}.$
\end{cor}

The above results are in complete agreement with the recursive  relations
between the wave functions and the explicit solution in terms of
the Mellin-Barnes representation obtained previously in \cite{KL2} by the
QISM-based approach.

\subsection{The hyperbolic Sutherland model}

In this section we give a short derivation of the integral representation
 for a  wave function of the hyperbolic
Sutherland model.

The quantum hyperbolic Sutherland model corresponding
to $\frak{gl}(N)$  is a integrable
quantum mechanical system with $N$ degrees of freedom (see \cite{OP2} and
references therein).
It has $N$ mutually commuting Hamiltonians, the first two of which
are given by
\be\label{cm00}
h_1 =\sum_{n=1}^Np_n\,,\\
h_2 = \sum_{m<n}\Big\{p_np_m+\frac{\hbar^2/4}
{\sinh^2(x_m-x_n)}\Big\}.
\ee
The eigenfunction  of the
Hamiltonians of the quantum hyperbolic Sutherland model
is defined by the
equations:
\be\label{cm0}
\sum_{n=1}^Np_n\,\Psi=\sigma_1(\bgamma_\N)\Psi\,,\\
\sum_{m<n}\Big\{p_np_m+\frac{\hbar^2/4}
{\sinh^2(x_m-x_n)}\,\Big\}\Psi=\sigma_2(\bgamma_{\N})\Psi\,.
\ee

In a similar fashion to the open Toda chain, the  hyperbolic
Sutherland model can
be explicitly solved using the Representation Theory approach.
Let $V$ be a representation of ${\cal U}(\frak{gl}(N,\RR))$ and
$v_\N$ be a vector in $V$ such that
 \be\label{cm1}
\hspace{2cm}
(E_{n,n+1}-E_{n+1,n})v_{\N}=0\,,\hspace{1cm}(n=1,\ldots, N-1).
\hspace{-2cm}
\ee
Consider the matrix element  $\Phi(x_1,\ldots,x_\N)$
of the form
\be\label{cm3}
\Phi=\<v_{\N},e^{-\sum\limits_{k=1}^Nx_kE_{kk}}v_{\N}\>\,.
\ee
One can show that (\ref{cm3}) gives (up to a simple factor)
the eigenfunction of the Hamiltonian of the hyperbolic
Sutherland model. To derive this,   consider the more general
matrix elements:
\be\label{cm3'}
\Phi_r=\<v_{\N},e^{-\sum\limits_{k=1}^Nx_kE_{kk}}a_{\N r}
v_{\N}\>\,,
\ee
where $a_{\N r}$ are the coefficients of the expansion
of $\ca_\N(\la)$
(see (\ref{cas1})): ${\cal A}_{\N}(\la)=
\la^N-\la^{N-1}a_{\N1}+
\la^{N-2}a_{\N2}\,+\ldots$. The first two have the form:
\be\label{cm4}
a_{\N 1}=i\hbar\sum_{n=1}^NE_{nn}\,,\\
a_{\N 2}=-\hbar^2\Big\{\sum_{m<n}E_{mm}E_{nn}-
\sum_{n=1}^N\rho^{(\N)}_nE_{nn}
-\sum_{m<n}E_{nm}E_{mn}+\sigma_2(\brho^{(\N)})\Big\}.
\ee

This should be compared with  the discussion of the
open Toda chain (\ref{pair1}).
The same line of reasoning shows that the matrix element
(\ref{cm3}) satisfies the differential equations
\be\label{cm10}
-i\hbar\sum_{n=1}^N\d_{x_n}\Phi=\sigma_1(\bgamma_\N)\Phi\,,\\
\hspace{-0.6cm}
-\hbar^2\sum_{m<n}\Big\{\d_{x_m}\d_{x_n}-{\textstyle\frac{1}{2}}
\coth(x_m-x_n)(\d_{x_m}-\d_{x_n})\Big\}\Phi=
\Big(\sigma_2(\bgamma_{\N})+
\hbar^2\sigma_2(\brho^{(\N)})\Big)\Phi\,.
\hspace{-1cm}
\ee
Actually, $\Phi$ is the $GL(N,\RR)$
zonal spherical function \cite{Gelf}-\cite{HC}.

Finally, the wave function $\Psi$ defined by
\be\label{cm11}
\Psi=\Phi\prod_{j<k}\sinh^{1/2}(x_j-x_k)
\ee
satisfies the equations (\ref{cm0}).

Now, we give a new explicit  integral formula for  the eigenfunction
of the hyperbolic Sutherland model using the
GZ representation discussed in section 2.
First, let us find the vector $v_\N$ satisfying (\ref{cm1}).
\begin{lem}
Let
\be\label{cm13}
\phi_n(\bgamma_n,\bgamma_{n+1})=\prod_{k=1}^{n}\prod_{m=1}^{n+1}
\hbar^{\frac{\gamma_{nk}-\gamma_{n+1,m}}{2i\hbar}+\frac{1}{4}}\;
\Gamma\Big(\frac{\gamma_{nk}-\gamma_{n+1,m}}{2i\hbar}+
\frac{1}{4}\Big).
\ee
Then the vector
\be\label{cm14}
v_{\N}=
e^{-\frac{\pi}{2\hbar}\sum\limits_{n=1}^{N-1}(n-1)
\sum\limits_{j=1}^n\gamma_{nj}}
\prod_{n=1}^{N-1}\phi_n(\bgamma_n,\bgamma_{n+1})
\ee
satisfies the  equations (\ref{cm1}).
\end{lem}
{\bf Proof.} Straightforward check.\square

\smallskip\noindent
Now we have the following
\begin{te}
The solution to (\ref{cm0}) admits the integral representation:
\be\label{su-f}
\Psi_{\gamma_{\N1},\ldots,\gamma_{\N\N}}(x_1,\ldots,x_\N)\,=
\prod_{j<k}\sinh^{1/2}(x_j-x_k)\\
\times \int\limits_{\RR^{\frac{N(N-1)}{2}}}
\prod_{n=1}^{N-1}
\frac{\prod\limits_{k=1}^n\prod\limits_{m=1}^{n+1}
\Big|\Gamma(\frac{\gamma_{nk}-\gamma_{n+1,m}}{2i\hbar}+
\frac{1}{4})\Big|^2}
{\prod\limits_{s<p}\left|
\Gamma(\frac{\gamma_{ns}-\gamma_{np}}{i\hbar})\right|^2}
\;e^{\frac{i}{\hbar}
\sum\limits_{n,j=1}^N(\gamma_{nj}-\gamma_{n-1,j})x_n}
\prod_{\stackrel{\scriptstyle n=1}{j\leq n}}^{N-1}d\gamma_{nj}\,.
\hspace{-1.5cm}
\ee
\end{te}
{\bf Proof.} The only non-obvious thing is the convergence of the
integral. However, the integral converges absolutely
due to the inequality (\ref{absol}).\square

\smallskip\noindent
{\bf Example.}
Let $N=2$. Then
\be\label{e1}
\Psi_{\gamma_{21},\gamma_{22}}(x_1,x_2)\,=\sinh^{1/2}(x_1-x_2)\\
\times\; e^{\frac{i}{\hbar}(\gamma_{21}+\gamma_{22})x_2}
\int\limits_{\RR}
\prod\limits_{m=1}^2
\Big|\Gamma\Big(\frac{\gamma_{11}-\gamma_{2m}}{2i\hbar}+
\frac{1}{4}\Big)\Big|^2
\;e^{\frac{i}{\hbar}\gamma_{11}(x_1-x_2)}d\gamma_{11}.
\ee
Assuming $x_1\geq x_2$, one can close the contour in the
upper half-plane thus calculating the integral by residues.
The answer is given in terms
of hypergeometric functions as follows:
\be
\Psi_{\gamma_{21},\gamma_{22}}(x_1,x_2)=\sinh^{1/2}(x_1-x_2)
\frac{4\pi^{\frac{1}{2}}\hbar}
{\cosh\frac{\pi}{2\hbar}(\gamma_{21}-\gamma_{22})}\,
e^{\frac{i}{\hbar}(\gamma_{21}+\gamma_{22})(x_1+x_2)}\\
\hspace{-2cm}
\times\; e^{-\frac{1}{2}(x_1-x_2)}\left\{
e^{\frac{i}{2\hbar}(\gamma_{21}-\gamma_{22})(x_1-x_2)}
\frac{\Gamma(-\frac{\gamma_{21}-\gamma_{22}}{2i\hbar})}
{\Gamma(\frac{1}{2}-\frac{\gamma_{21}-\gamma_{22}}{2i\hbar})}
F(\ts{1}{2},\ts{1}{2}+\ts{\gamma_{21}-\gamma_{22}}{2i\hbar},
1+\ts{\gamma_{21}-\gamma_{22}}{2i\hbar};e^{-2(x_1-x_2)})\right.
\hspace{-2cm}\\
\left.
\hspace{-1cm}
+\;e^{-\frac{i}{2\hbar}(\gamma_{21}-\gamma_{22})(x_1-x_2)}
\frac{\Gamma(\frac{\gamma_{21}-\gamma_{22}}{2i\hbar})}
{\Gamma(\frac{1}{2}+\frac{\gamma_{21}-\gamma_{22}}{2i\hbar})}
F(\ts{1}{2},\ts{1}{2}-\ts{\gamma_{21}-\gamma_{22}}{2i\hbar},
1-\ts{\gamma_{21}-\gamma_{22}}{2i\hbar};e^{-2(x_1-x_2)})
\right\}.
\ee
Up to a trivial factor this solution coincides with the Legendre
function $P_{\frac{\gamma_{21}-\gamma_{22}}{2i\hbar}-
\frac{1}{2}}(\cosh(x_1\!-\!x_2))$
(see \cite{BE}, formulas 3.2(9), 3.2(27)), namely,
\be
\Psi_{\gamma_{21},\gamma_{22}}(x_1,x_2)\,=\sinh^{1/2}(x_1-x_2)\\
\times\;\frac{4\pi\hbar}
{\cosh\frac{\pi}{2\hbar}(\gamma_{21}-\gamma_{22})}\,
e^{\frac{i}{\hbar}(\gamma_{21}+\gamma_{22})(x_1+x_2)}
P_{\frac{\gamma_{21}-\gamma_{22}}{2i\hbar}-\frac{1}{2}}
(\cosh(x_1-x_2))\,.
\ee
Note,  in the limit $x_1- x_2 \rightarrow\infty$
the wave function (\ref{e1}) has the asymptotic behaviour
\be
\Psi_{\gamma_{21},\gamma_{22}}(x_1,x_2)\sim
c(\gamma_{21},\gamma_{22})
e^{\frac{i}{\hbar}(\gamma_{21}x_1+\gamma_{22}x_2)}
+c(\gamma_{22},\gamma_{21})
e^{\frac{i}{\hbar}(\gamma_{22}x_1+\gamma_{21}x_2)}\,,
\ee
where
\be
c(\la_1,\la_2)=\frac{\Gamma(-\frac{\la_1-\la_2}{2i\hbar})}
{\Gamma(\frac{1}{2}-\frac{\la_1-\la_2}{2i\hbar})}
\ee
is the Harish-Chandra function.
\begin{rem}
Since the integrand in (\ref{su-f}) is written in the factorized form,
it is clear that there exists a recursive relation, similar to (\ref{rec}),
between the $n$ and $n-1$ particle eigenfunctions of the
Sutherland model. This fact is the direct consequence of the
Gelfand-Zetlin inductive procedure.
\end{rem}

\section{Connection with the QISM-based approach}

In section 5  we have outlined the derivation of the integral
representation of the wave function of an  open Toda chain based on
the GZ representation of ${\cal U}(\frak{gl}(N,\RR))$.
As was stressed in section 4, this representation bears a direct
relationship to $Y(\frak{gl}(N,\RR))$. This provides us with
another point of view on the way to obtain the wave functions
for the open Toda chain. In particular, the algebra generated by
the images of the polynomials $\ca_n(\la)$, $\cb_n(\la)$,
$\cc_n(\la)$ under the homomorphism
$\pi_\N: Y(\frak{gl}(N))\rightarrow{\cal U}(\frak{gl}(N))$
can successfully be used to derive the recursive relations for
the wave functions.

We begin with a short summary of the QISM-based approach to the
$GL(N,\RR)$ Toda chain due to Sklyanin \cite{Skl2}. The model is
described with the help of the auxiliary $2\times 2$ matrix
Lax operator with spectral parameter:
\be\label{ln}
L_n(\la)\,=\,
\left(\begin{array}{cc}\la-p_n & e^{-x_n}\\ -e^{x_n} & 0
\end{array}\right),
\ee
where $p_n=-i\hbar\d_{x_n}$.

The operator (\ref{ln}) satisfies the usual quadratic relation
\be\label{rll}
R_{2}(\la-\mu)(L_n(\la)\otimes I)(I\otimes L_n(\mu))=
(I\otimes L_n(\mu))(L_n(\la)\otimes I)R_{2}(\la-\mu)
\ee
with  the rational $4\times 4$ $R$-matrix (see (\ref{r-mat})).
The monodromy matrix
\be\label{mon1}
T_n(\la)=L_{n}(\la)\ldots
L_1(\la)\,\equiv\, \left(\begin{array}{cc}
A_{n}(\la) & B_{n}(\la)\\
C_{n}(\la) & D_{n}(\la)\end{array}\right)
\ee
also satisfies the  relation:
\be\label{rtt}
R_{2}(\la-\mu)(T_{n}(\la)\otimes I)(I\otimes T_{n}(\mu))=
(I\otimes T_{n}(\mu))(T_{n}(\la)\otimes I)R_{2}(\la-\mu).
\ee
In particular, the following commutation relations hold:
\be\label{carel}
[A_{n}(\la),A_{n}(\mu)]=[B_{n}(\la),B_{n}(\mu)]=
[C_{n}(\la),C_{n}(\mu)]=0\,,\\
(\la-\mu+i\hbar)A_{n}(\la)B_{n}(\mu)\,=\,
(\la-\mu)B_{n}(\mu)A_{n}(\la)+i\hbar A_{n}(\mu)B_{n}(\la)\,,\\
(\la-\mu+i\hbar)A_{n}(\mu)C_{n}(\la)\,=\,
(\la-\mu)C_{n}(\la)A_{n}(\mu)+i\hbar A_{n}(\la)C_{n}(\mu)\,.
\ee
The polynomial $\!A_n(\la)\!$ is a generating function for
the quantum Hamiltonians of the $\!GL(n,\RR)\!$  Toda chain, and
the wave functions $\psi_{\bgamma_n}(x_1,\ldots,x_n),\;n=1,\ldots,N$
are defined as solutions of the differential equations
\be\label{km4}
A_n(\la)\psi_{\gamma_{n1},\ldots,\gamma_{nn}}=
\prod_{j=1}^n(\la-\gamma_{nj})\;
\psi_{\gamma_{n1},\ldots,\gamma_{nn}}\,.
\ee
The idea of the QISM approach is to solve iteratively these equations
for  $\psi_{\bgamma_n},\, n=1,\ldots,N$ using the additional operators
$B_n(\la)$ and $C_n(\la)$. From the results of \cite{KL2} (eq. (2.7)),
the following expressions for the operators
$A_n(\la),\,B_n(\la),$ and $C_n(\la)$ can be obtained:
\be\label{km12}
A_n(\la)=\prod_{j=1}^n(\la-\gamma_{nj})\,,\\
B_n(\la)=i^{1+n}
\sum_{j=1}^n\prod_{s\neq j}\frac{\la-\gamma_{ns}}
{\gamma_{nj}-\gamma_{ns}}\;e^{-i\hbar\d_{\gamma_{nj}}}\,,\\
C_n(\la)=i^{1-n}
\sum_{j=1}^n\prod_{s\neq j}\frac{\la-\gamma_{ns}}
{\gamma_{nj}-\gamma_{ns}}\;e^{i\hbar\d_{\gamma_{nj}}}\,.
\ee
A simple check shows that, thus defined, the operators
(\ref{km12}) satisfy
the relations (\ref{carel}). Hence, one establishes
the connection between
the generators ${\cal A}_n(\la),\,{\cal B}_n(\la)$
and ${\cal C}_n(\la)$
entering the description of the GZ representation (\ref{bza}),
and the corresponding QISM generators (\ref{km12}):
\be
A_n(\la)={\cal A}_n (\la)\,,\\
B_n(\la)=s^{-1}_n\circ {\cal B}_n (\la)\circ s_n\,,\\
C_n(\la)=s^{-1}_{n-1}\circ{\cal C}_n(\la)\circ s_{n-1}\,,
\ee
where $s_n$ is defined by (\ref{wv'}). Therefore, the expressions
(\ref{carel}) should be compared with (\ref{cw1}).

At the end of this section, we illustrate the connection between
the differential operator $A_n(\la)$ entering in (\ref{mon1}) and
the generator ${\cal A}_n(\la)$ defined by (\ref{cas1}).
First, note that (\ref{ln}), (\ref{mon1}) define the relations between
the operators $A_n(\la)$ for different levels:
\be\label{rec3}
A_n(\la)=(\la\!-\!p_{n})A_{n-1}(\la)-
e^{x_{n-1}-x_n}A_{n-2}(\la)\,,
\ee
where $n=1,\ldots,N$ and $A_{-1}=0,\,A_0=1$.
\begin{prop}
Let $\psi_{\gamma_{\N1},\ldots,\gamma_{\N\N}}$ be defined by
(\ref{pair}). Then the operators $A_n(\la)$ and $\ca_n(\la)$ are related
by
\be\label{h2}
A_n(\la)\psi_{\gamma_{\N1},\ldots,\gamma_{\N\N}}
=e^{-{\bfit x}\cdot\brho^{(N)}}
\langle {w}'_{\N},e^{-\sum\limits_{k=1}^Nx_kE_{kk}}
{\cal A}_n(\la-{\textstyle\frac{i(N-n)\hbar}{2}})
w_{\N} \rangle.
\ee
\end{prop}
{\bf Proof.} We prove the statement by showing that the operators
$A_n(\la)$ defined by (\ref{h2}) satisfy the recursive relations
(\ref{rec3}). Consider the relation (\ref{char1}). The generators
$E_{nn}$ and $E_{n,n-k+1}$ commute with ${\cal A}_{n-1}(\la)$ and
${\cal A}_{n-k}(\la)$, respectively. Therefore, due to
(\ref{inv5}) and (\ref{fww'}),
\be\label{h4}
\langle w'_{\N},e^{-\sum\limits_{k=1}^Nx_kE_{kk}}
{\cal A}_n({\textstyle\la-\frac{i(N-n)\hbar}{2}})w_{\N} \rangle \\
=\;(\la-i\hbar\rho^{(\N)}_n+i\hbar\d_{x_n})
\langle  w'_{\N},e^{-\sum\limits_{k=1}^Nx_kE_{kk}}
{\cal A}_{n-1}(\la-{\textstyle\frac{i(N-n+1)\hbar}{2}})
w_{\N} \rangle  \\
-\;e^{x_{n-1}-x_n} \langle
w'_{\N},e^{-\sum\limits_{k=1}^Nx_kE_{kk}}{\cal A}_{n-2}
(\la-{\textstyle\frac{i(N-n+2)\hbar}{2}})w_{\N} \rangle.
\ee
\square

We leave a more systematic discussion of the interrelation
between the two sets of operators $\ca_n(\la)$,$\cb_n(\la)$,$\cc_n(\la)$
and $A_n(\la)$,$B_n(\la)$,$C_n(\la)$ for a future publication.

\end{document}